\documentclass[12pt]{article}

\usepackage{graphicx}  
\usepackage{amsmath}  
\usepackage{amssymb,amsthm}
\usepackage{color}
\usepackage{tikz}
\usepackage{pgfplots}

\usepackage{csquotes}
\usepackage[numbers]{natbib}
\usepackage{notes2bib}

\bibnotesetup{
	note-name = ,
	use-sort-key = false
}

\newcommand{\R}{{\mathbb R}}\newcommand{\N}{{\mathbb N}}

\newcommand{\C}{{\mathbb C}}

\let\epsilon\varepsilon
\let\hat\widehat

\newtheorem{theorem}{Theorem}[section]\newtheorem{lemma}[theorem]{Lemma}
\newtheorem{definition}[theorem]{Definition}

\newtheorem{remark}[theorem]{Remark}\newtheorem{example}[theorem]{Example}

\title{Global existence  for long wave  Hopf unstable  \\ spatially extended systems  \\ with 
a conservation law}
\author{Nicole Gauss$^1$, Anna Logioti$^1$, Guido Schneider$^1$, \\[2mm] and Dominik Zimmermann$^1$\\
{\small
$^1$Institut f\"ur Analysis, Dynamik und Modellierung,} \\ {\small Universit\"at Stuttgart, Pfaffenwaldring 57, 
70569 Stuttgart, Germany}
\\ {\small email: firstname.lastname@mathematik.uni-stuttgart.de}}

\usepackage{filecontents}

\begin{document}

\maketitle

\begin{abstract}
We are interested in  reaction-diffusion systems, with a  conservation law,
exhibiting a Hopf bifurcation at the spatial wave number $ k = 0 $. 
With the help of a multiple scaling perturbation ansatz a Ginzburg-Landau equation 
coupled to a scalar   conservation law can be derived as an amplitude system for the 
approximate 
description of  the dynamics of the original   reaction-diffusion  system near the first instability.
We use the amplitude system to show  the global existence of all  solutions 
starting in a small neighborhood of the weakly unstable ground state
for original systems  posed 
on a large spatial interval  with periodic boundary conditions.
\end{abstract}


%

\section{Introduction}

We consider reaction-diffusion systems for $ u $ 
with $ u(x,t) \in \R^d $ for $ d \geq 2 $ coupled to a diffusive conservation law 
for $ v $  with $ v(x,t) \in \R $, namely 
\begin{eqnarray} \label{rd1}
\partial_t u & = & D \partial_x^2 u + f(u,v) , \\
\partial_t v & = & d_v \partial_x^2 v + \partial_x^2 g(u) ,
\label{rd2}
\end{eqnarray}
where $ x \in \R $, $ t \geq 0 $, $ D $ a diagonal diffusion matrix with entries $ d_j > 0 $ for 
$ j = 1,\ldots,d $,  $ d_v > 0 $ a scalar diffusion coefficient,  and $ f: \R^d \times \R \to \R^d $
and $ g: \R^d  \to \R $  smooth reaction terms  with 
$$
f(u,v) = \mathcal{O}(|u|(1+|u|+|v|)   \qquad  \textrm{and}  \qquad 
g(u) =   \mathcal{O}(|u|^2)  $$ such that 
$ (u,v)  = (0,v^*) $ is a stationary solution for any constant  $ v^* \in \R $.
As a consequence of the conservation law form
the  spatial integral of $ v $ is conserved in time. 
The fact that $ g $ only depends on $ u $ or that $ D $ is a diagonal matrix are no 
restrictions w.r.t. our purposes. For a detailed discussion about that see Section \ref{secdisc}.

We are interested in the behavior of \eqref{rd1}-\eqref{rd2}
close to the stationary solutions, w.l.o.g. for our purposes take  $ (u,v)  = (0,0) $.
The linearization of   \eqref{rd1}-\eqref{rd2} at $ (0,0) $,
\begin{eqnarray} \label{lrd1}
\partial_t u & = & \Lambda_ u  u = D \partial_x^2 u + \partial_u f(0,0) u  ,\\
\partial_t v & = & \Lambda_ v v= d_v \partial_x^2 v , 
\label{lrd2}
\end{eqnarray}
is solved by $ u(x,t) = e^{ikx + \lambda t} \widehat{u} $ and 
$ v(x,t) = e^{ikx + \lambda t} \widehat{v} $ where $ \lambda \in \C $, 
$ \widehat{u}  \in \C^d $, and $ \widehat{v}  \in \C $ are determined by 
\begin{eqnarray} \label{lrd1ev}
 \lambda \widehat{u} & = & - D k^2 \widehat{u} + \partial_u f(0,0) \widehat{u} , \\
 \lambda \widehat{v} & = & - d_v k^2 \widehat{v} . 
\label{lrd2ev}
\end{eqnarray}
We find $ d $ curves of eigenvalues $ \lambda_j =  \lambda_j(k)  $ ordered  
as $ \textrm{Re} \lambda_1(k) \geq \ldots \geq \textrm{Re} \lambda_d(k) $ for 
\eqref{lrd1ev} and $ \lambda_0(k) = - d_v k^2 $ for \eqref{lrd2ev}.
The associated normalized eigenvectors or normalized generalized eigenvectors are denoted by $  \widehat{U}_j \in \mathbb{C}^d$
for $ j = 0,\ldots, d $.

We assume that  \eqref{rd1}-\eqref{rd2} depends on a parameter $ \widetilde{\alpha} $
and that for $ \widetilde{\alpha} = \widetilde{\alpha}_c $ we have the following spectral situation.
\medskip

{(\bf
 Spec)}  There is an $ \omega_0 > 0 $ such that 
$ \textrm{Re}  \lambda_j(0)|_{\widetilde{\alpha} = \widetilde{\alpha}_c} = \lambda'_j(0)|_{\widetilde{\alpha} = \widetilde{\alpha}_c} = 0 $, 
$ \textrm{Re}  \lambda''_j(0)|_{\widetilde{\alpha} = \widetilde{\alpha}_c} < 0 $  for $ j = 1,2 $ and $ \textrm{Im}  \lambda_1(0)|_{\widetilde{\alpha} = \widetilde{\alpha}_c} = 
- \textrm{Im}  \lambda_2(0)|_{\widetilde{\alpha} = \widetilde{\alpha}_c} = 
\omega_ 0 $. Moreover, all other eigenvalues $ \lambda_j|_{\widetilde{\alpha} = \widetilde{\alpha}_c} $ 
for $ j =1 ,\ldots , d $
have a  negative real part.
Finally, we assume that $ \partial_{\widetilde{\alpha}}  \textrm{Re}  \lambda_1(0)|_{\widetilde{\alpha} = \widetilde{\alpha}_c} > 0 $.
\medskip

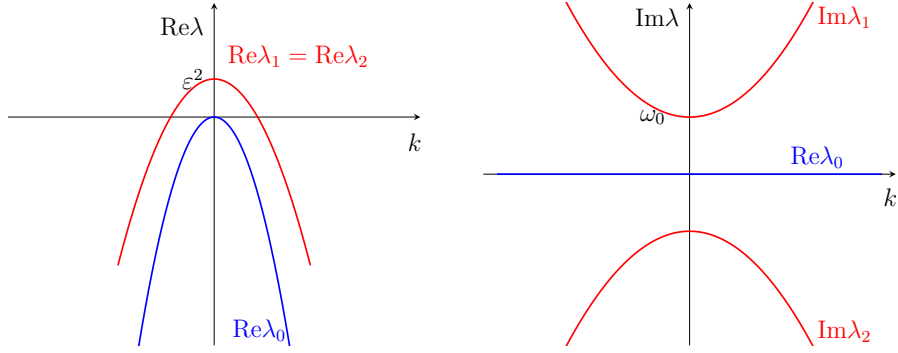
\begin{figure}[h] 
   \centering
     \begin{tikzpicture}[rotate=0,scale=0.8]
    \begin{axis}[
	xmin=-1.5, xmax=1.5,
	ymin=-0.6, ymax=0.3,
	axis lines=center,
	ticks=none,
    ]
      \addplot+[blue,thick, domain=-0.6:0.6,samples=100,no marks,opacity=0.5] {-2*x^2};
         \addplot+[red,thick, domain=-0.7:0.7,samples=100,no marks,opacity=0.5] {0.1-x^2};
           \node[below,black] at (axis cs:1.46,-0.02) {$ k $};
      \node[] at (axis cs:-0.35,0.45) {$\textrm{Re}(\lambda_+)$};
      \node[left=4pt] at (axis cs:0.05,0.24) {$ \textrm{Re} \lambda $};
         \node[red, left=4pt] at (axis cs:1.25,0.16) {$ \textrm{Re} \lambda_1 = \textrm{Re} \lambda_2$};
           \node[blue, left=4pt] at (axis cs:0.65,-0.56) {$ \textrm{Re} \lambda_0 $};
        \node[left=4pt] at (axis cs:0.05,0.1) {$ \varepsilon^2 $};
    \end{axis}
  \end{tikzpicture}\qquad
    \begin{tikzpicture}[rotate=0,scale=0.8]
    \begin{axis}[
	xmin=-1.5, xmax=1.5,
	ymin=-0.6, ymax=0.6,
	axis lines=center,
	ticks=none,
    ]
      \addplot+[blue,thick, domain=-1.4:1.4,samples=100,no marks,opacity=0.5] {0};
         \addplot+[red,thick, domain=-1.4:1.4,samples=100,no marks,opacity=0.5] {-0.2-x^2/2};
           \addplot+[red,thick, domain=-1.4:1.4,samples=100,no marks,opacity=0.5] {0.2+x^2/2};
           \node[below,black] at (axis cs:1.46,-0.02) {$ k $};
      \node[left=4pt] at (axis cs:0.05,0.55) {$ \textrm{Im} \lambda $};
        \node[left=4pt] at (axis cs:-0.05,0.2) {$ \omega_0 $};
           \node[red, left=4pt] at (axis cs:1.45,0.55) {$ \textrm{Im} \lambda_1 $};
                \node[red, left=4pt] at (axis cs:1.45,-0.55) {$ \textrm{Im} \lambda_2 $};
           \node[blue, left=4pt] at (axis cs:1.25,0.06) {$ \textrm{Re} \lambda_0 $};
    \end{axis}
  \end{tikzpicture}\qquad

   \caption{The relevant spectral curves  of the linearization around  the  trivial solution
   plotted as a function over the Fourier wave numbers for $ \widetilde{\alpha} - \widetilde{\alpha}_c  = \varepsilon^2 > 0 $. The left panel shows the real part of the eigenvalue curves $ \lambda_0 $ (in blue), $ \lambda_1 $, 
   and $ \lambda_2 $ (both in red), the right panel shows the imaginary  part.
}
   \label{figspec}
\end{figure}

 For \eqref{rd1}-\eqref{rd2} from the assumption  {(\bf
 Spec)}   a spectral situation follows as sketched in 
 Figure \ref{figspec}.

 {\bf Notation.} In order to make the notation more intuitive in the following we use the index $ -1 $ instead of $ 2 $, i.e., for example we write $ \lambda_{-1} = \lambda_2 $.

 We introduce the bifurcation parameter 
$\varepsilon^2 = \widetilde{\alpha} - \widetilde{\alpha}_c  $ and  insert the ansatz 
\begin{eqnarray} \label{glans1}
u(x,t) & = & \varepsilon A_1(X, T ) e^{i \omega_0 t} \widehat{U}_1(0) + \textrm{c.c.}
+ \mathcal{O}(\varepsilon^2) , \\
v(x,t) & = & \varepsilon^2 B_0(X, T ) ,\label{glans2}
\end{eqnarray}
with $ X = \varepsilon  x $, $ T = \varepsilon^2 t $, $  B_0(X, T ) \in \mathbb{R} $,
and $  A_1(X, T ) \in \mathbb{C} $ in  \eqref{rd1}-\eqref{rd2}.
We obtain the system of amplitude equations 
\begin{eqnarray} \label{as1intro}
\partial_T A_1  & = & a_0 \partial_X^2 A_1 + a_1 A_1+  a_2  A_1 B_0 -  a_3  A_1 |A_1|^2,\\
\partial_T B_0& =  &  b_0 \partial_X^2 B_0 +  b_1 \partial_X^2 (|A_1|^2), \label{as2intro}
\end{eqnarray}
with coefficients $ a_0, a_3 \in \mathbb{C} $, $ a_1,a_2,b_0,b_2 \in \mathbb{R} $, satisfying $ \textrm{Re} a_0 > 0 $,  $ b_0 > 0 $, $ a_1 > 0 $, and
$  \textrm{Re} a_3 > 0 $,
consisting of a Ginzburg-Landau equation for $ A_1 $ coupled to a
scalar conservation law for $ B_0 $.
The amplitude function $ A_1 $ describes the oscillatory modes concentrated a $ k = 0 $ and 
$ B_0 $  the conservation law modes concentrated at $ k = 0 $.

\begin{example}\label{exam11}{\rm 
In order to make this introduction less abstract 
the derivation of the amplitude system 
will be explained 
for the following toy problem 
\begin{eqnarray*}
\partial_t u_1 & = & \partial_x^2 u_1 + i \omega_0 u_1 + \varepsilon^2 u_1 
+ u_1^2 + u_1 u_{-1} + u_{-1}^2 + v u_1 + v u_{-1} - u_1^2 u_{-1}, \\
\partial_t u_{-1}  & = & \partial_x^2 u_{-1}  - i \omega_0 u_{-1}  + \varepsilon^2 u_{-1}  
+ u_1^2 + u_1 u_{-1} + u_{-1}^2 + v u_1 + v u_{-1}- u_{-1}^2 u_{1}, \\
\partial_t v & = & \partial_x^2 v +   \partial_x^2 ( u_1 u_{-1}) ,
\end{eqnarray*}
with $ u_{-1} = \overline{u_1} $. 
Although it is not of the form of \eqref{rd1}-\eqref{rd2}, it shares 
essential properties with \eqref{rd1}-\eqref{rd2}, in particular, it has 
qualitatively a spectral picture as plotted in Figure \ref{figspec}.
We make the ansatz 
\begin{eqnarray*}
u_1 (x,t) & = &  \varepsilon A_1(X, T ) e^{i \omega_0 t} 
+ \varepsilon^2 A_{1,0}(X, T ) \\ &&
+ \varepsilon^2 A_{1,2}(X, T ) e^{2 i \omega_0 t} + \varepsilon^2 A_{1,-2}(X, T ) e^{- 2 i \omega_0 t} , \\
u_{-1} (x,t) & = &  \varepsilon A_{-1}(X, T ) e^{-i \omega_0 t} + \varepsilon^2 A_{-1,0}(X, T )
\\ &&
+ \varepsilon^2 A_{-1,2}(X, T ) e^{2 i \omega_0 t} + \varepsilon^2 A_{-1,-2}(X, T ) e^{- 2 i \omega_0 t} , \\
v(x,t) & = & \varepsilon^2 B_0(X, T ) ,
\end{eqnarray*}
with $ A_{-1} = \overline{A_1} $, etc.
For the $ u_1 $-equation we find:
\begin{eqnarray*}
\varepsilon^3 e^{i\omega_0 t} & : & \partial_T A_1  =  \partial_X^2 A_1 + A_1 
+ B_{0} A_1  
\\ && \qquad \qquad + 2 A_{1,0} A_1  + A_{1,2} A_{-1}
+ A_{-1,0} A_1 + 2 A_{-1,2} A_{-1}
- A_1^2 A_{-1} , \\
\varepsilon^2 e^{2 i\omega_0 t} & : & 2 i \omega_0 A_{1,2} = i \omega_0 A_{1,2} + A_1^2 , \\
\varepsilon^2 e^{0 i\omega_0 t} & : & 0= i \omega_0 A_{1,0} + A_1 A_{-1} , \\
\varepsilon^2 e^{- 2 i\omega_0 t} & : & - 2 i \omega_0 A_{1,-2} = i \omega_0 A_{1,-2} + A_{-1}^2 .
\end{eqnarray*}
For the $ u_{-1} $-equation we find similar equations and for the 
$ v $-equation we obtain:
\begin{eqnarray*}
\varepsilon^4 & : & \partial_T B_0  = \partial_X^2 B_{0} +  \partial_X^2(A_1 A_{-1}) .
\end{eqnarray*}
If we eliminate the $ A_{j,0} $ and $ A_{j,2} $ by the above algebraic equations we find 
\begin{eqnarray*}
 \partial_T A_1  & = &  \partial_X^2 A_1 + A_1 
+ B_{0} A_1  - \gamma_3 |A_1|^2 A_{1} ,\\
 \partial_T B_0  & = & \partial_X^2 B_{0} +  \partial_X^2(|A_1|^2) ,
\end{eqnarray*}
with 
\begin{equation} \label{gamma3} 
- \gamma_3 =  - \frac{2}{i \omega_0} + \frac{1}{i \omega_0} + \frac{1}{i \omega_0} + \frac{2}{3 i \omega_0} -1 =  -1 + \frac{2 }{3 i \omega_0} .
\end{equation}
\mbox{}\hfill$\lrcorner$
}\end{example}

In order to establish the global existence and uniqueness for \eqref{as1intro}-\eqref{as2intro}, in the following we assume 
\medskip 

{\bf (Coeff)}  The coefficients $ a_0, \ldots , b_1 $ of \eqref{as1intro}-\eqref{as2intro} satisfy 
 for the  normalized System \eqref{as5}-\eqref{as6},
subsequently computed in Remark \ref{rems}, that $ 1+  \alpha^{-1} \beta > 0 $.
\medskip 

Using the same multiple scaling analysis,  in \cite{Schn98JNS},
in case of no conservation law, i.e., in case $ v = 0 $ and without the $ v $-equation
in \eqref{rd1}-\eqref{rd2}, a Ginzburg-Landau equation was derived, and it was shown that all small solutions develop in such a way that they can be approximated after a certain time by the solutions of the Ginzburg-Landau equation. The proof differs essentially from the case when the bifurcating pattern is oscillatory in space
which is based on mode-filters and a detailed analysis of the mode interactions.
See \cite[\S 10]{SU2017book} for an overview.
In contrast the  proof of \cite{Schn98JNS} is based on normal form methods. As a consequence of the results of \cite{Schn98JNS}, the global existence in time of all small bifurcating solutions and the upper-semicontinuity of the rescaled original system attractor towards the associated Ginzburg-Landau attractor follows. The result of \cite{Schn98JNS}
applies for instance to the Brusselator, the Schnakenberg, the Gray-Scott or 
the Gierer-Meinhardt model, cf. \cite{SW22}.

It is the purpose of this paper to prove a similar global existence result for 
 \eqref{rd1}-\eqref{rd2}, i.e., in case of an additional conservation law,
 with the help of the   amplitude system \eqref{as1intro}-\eqref{as2intro}.
 
 This  question
turns out to be very challenging for the following reason.
Since $ (A_1,B_0) = (0,B^*) $, with constants $ B^* \in \R $, 
 is an unbounded family  of stationary solutions
 for \eqref{as1intro}-\eqref{as2intro},
this amplitude system 
does not possess an exponentially absorbing ball if posed on the real line,
in contrast to a single Ginzburg-Landau equation if $ \textrm{Re} a_3 > 0 $.
However, 
assuming {\bf (Coeff)} 
an exponentially attracting ball exists in case of periodic boundary conditions, say  
\begin{equation} \label{ABper}
A_1(X,T) = A_1(X+2 \pi,T) \qquad \textrm{and} \qquad 
B_0(X,T) = B_0(X+2 \pi,T).
\end{equation} 
Then we have the existence of an absorbing ball and the global existence and uniqueness of solutions.
\begin{theorem} \label{thGL}
Consider the amplitude system \eqref{as1intro}-\eqref{as2intro} with periodic
boundary conditions \eqref{ABper} and assume that the coefficients 
$ a_0 ,\ldots , b_1 $ satisfy the  condition {\bf (Coeff)}.
Then for all $ s \geq 0 $ there exists a $ C_R > 0 $ such that for all $ C_1 > 0 $ there exists a $ T_0 > 0 $ 
such that to a given initial condition $  (A_1(\cdot,0),B_0(\cdot,0)) \in H^{s+1} \times 
H^s $ with $ \| A_1(\cdot,0) \|_{H^{s+1}} +  \| B_0(\cdot,0) \|_{H^s} \leq C_1 $
there exists a unique global solution 
$
(A_1,B_0) \in C([0,\infty), H^{s+1} \times H^s )
$ 
such that additionally
$ \| A_1(\cdot,T) \|_{H^{s+1}} +  \| B_0(\cdot,T) \|_{H^s} \leq C_R $ for all $ T \geq T_0 $.   
\end{theorem}
\begin{remark}\label{rem13n}{\rm 
In case of periodic boundary conditions the Sobolev space $ H^s $ can be embedded in the space  $ H^s_{l,u} $ of uniformly local Sobolev functions
for  $ s \geq 0 $ and so  in case of periodic boundary conditions, the existence of an absorbing ball in
$ H^{s+1}_{l,u} \times H^s_{l,u} $ for $ (A_1,B_0) $ follows, too. 
For the definition of the space $ H^s_{l,u} $ see the notations 
on Page \pageref{notat}.
}
\end{remark}
As already said we are interested in a similar result for the original system 
\eqref{rd1}-\eqref{rd2} using the existence of an exponentially attracting 
absorbing ball for the amplitude system \eqref{as1intro}-\eqref{as2intro} and
the fact that all solutions of  \eqref{rd1}-\eqref{rd2} develop in such a way 
that after a certain time they can be approximated by the solutions
of the amplitude system \eqref{as1intro}-\eqref{as2intro}.

The $ 2 \pi $-spatially periodic boundary conditions for the amplitude system \eqref{as1intro}-\eqref{as2intro} correspond to $ 2 \pi/\varepsilon $-spatially periodic boundary conditions for the original system 
\eqref{rd1}-\eqref{rd2}, 
i.e.,
\begin{equation} \label{uvper}
u(x,t) = u(x+2 \pi/\varepsilon,t) \qquad \textrm{and} \qquad 
v(x,t) = v(x+2 \pi/\varepsilon,t) .
\end{equation} 
Then for these periodic boundary conditions and small $ \varepsilon > 0 $ we have the  global existence and uniqueness of solutions for the original system 
\eqref{rd1}-\eqref{rd2}.

\begin{theorem} \label{thrd}
Consider the original  system \eqref{rd1}-\eqref{rd2} with periodic
boundary conditions \eqref{uvper} and assume that the coefficients 
$ a_0 ,\ldots , b_1 $ of the associated amplitude system \eqref{as1intro}-\eqref{as2intro}
satisfy the condition {\bf (Coeff)}.
Then for all $ n \geq 0 $ there exists a $ C_R > 0 $ and an $ \varepsilon_0 > 0 $  such that for all $ C_1 > 0 $ 
and all $ \varepsilon \in (0,\varepsilon_0) $,
there exists a $ t_0 = \mathcal{O}(1/\varepsilon^2)> 0 $ 
such that to a given initial condition $  (u(\cdot,0),v(\cdot,0)) \in H^{n+1} _{l,u} \times 
H^n_{l,u} $ with $ \| u(\cdot,0) \|_{H^{n+1}_{l,u}} + \varepsilon^{-1} \| v(\cdot,0) \|_{H^n_{l,u}} \leq C_1 \varepsilon$
there exists a unique global solution 
$
(u,v) \in C([0,\infty),H^{n+1} _{l,u} \times 
H^n_{l,u} )
$ 
such that additionally
$ \| u(\cdot,t) \|_{H^{n+1} _{l,u}} +  \varepsilon^{-1}  \| v(\cdot,t) \|_{H^n_{l,u}  } \leq C_R \varepsilon $ for all $ t \geq t_0 $.   
\end{theorem}

\begin{remark}{\rm 
Hence,  the global existence question can be answered positively at least 
 for original systems with periodic boundary 
 conditions  \eqref{uvper} which correspond in the amplitude system
 \eqref{as1intro}-\eqref{as2intro} to periodic boundary 
 conditions \eqref{ABper}. Since the $ L^2 $-norm of $ u = 1 $ on the interval $ [-\pi/\varepsilon,\pi/\varepsilon] $
grows as $ \mathcal{O}(1/\sqrt{\varepsilon}) $ with $ \varepsilon \to 0 $, 
Sobolev spaces are not adequate for controlling the norm and so
spaces have to be used where functions such as $ u= 1 $ 
can be  bounded independently of the small 
perturbation parameter $ 0 < \varepsilon \ll 1 $.
}
\end{remark}
\begin{remark}{\rm \label{rem15}
The three main ingredients of  the global existence proof
are {\bf (GL)}: the existence of an exponentially attracting absorbing ball of the amplitude system,
{\bf (APP)}:
an approximation result 
which shows that solutions of the original system \eqref{rd1}-\eqref{rd2} 
can be approximated on the natural $ \mathcal{O}(1/\varepsilon^2) $-time scale 
of  \eqref{as1intro}-\eqref{as2intro} of the amplitude system via the solutions of \eqref{as1intro}-\eqref{as2intro},
and {\bf (ATT)}:
an  attractivity result, which shows that solutions of \eqref{rd1}-\eqref{rd2} 
to initial conditions of order $ \mathcal{O}(\varepsilon) $ develop 
in such a way that after an $ \mathcal{O}(1/\varepsilon^2) $-time scale 
they are of a form which allows us to approximate them 
afterwards  by the solutions of \eqref{as1intro}-\eqref{as2intro}.
}
\end{remark}
\begin{remark}{\rm 
Approximation and attractivity results have been established  in \cite{HSZ11,SZ13,DSSZ16} in case of a Turing pattern forming systems coupled to a conservation law. Attractivity and approximation results 
in case of a simultaneous Turing and  a long wave Hopf bifurcation
 can be found in \cite{SW22}.}
 \end{remark}
 \begin{remark}{\rm
 The idea is as follows.  
   A neighborhood  of the origin of the 
   pattern forming system is mapped by the attractivity  {\bf (ATT)} into a set which can be 
   described by the amplitude system. The amplitude system 
   possesses  an exponentially attracting absorbing ball {\bf (GL)}.
   Therefore, by the approximation property {\bf (APP)} the original neighborhood  of the 
   pattern forming system is mapped after a certain time into itself.
   These a priori estimates combined with the local existence and uniqueness gives 
   the global existence and uniqueness of solutions of the pattern forming system in 
   a neighborhood of the weakly unstable origin.
}
 \end{remark}
 
 \begin{remark}{\rm 
Examples of reaction-diffusion systems \eqref{rd1}-\eqref{rd2}, falling into the class of systems we are interested in, are for instance the Brusselator, the Schnakenberg, 
the Gray-Scott  and the Gierer-Meinhardt model coupled to a conservation law 
coming for instance  from ecology. As an example we consider the Brusselator. 
The system, with the spatially homogeneous trivial equilibrium as origin, is given
by 
\begin{eqnarray}
\label{brussel2a}
\partial_t u_1 & = & d_1 \partial_x^2 u_1+ (b-1) u_1 + a^2 u_2  + f(u_1,u_2) ,\\
\partial_t u_2 & = & d_2 \partial_x^2 u_2 -b u_1 - a^2 u_2 
- f(u_1,u_2) , \label{brussel2b}
\end{eqnarray}
with nonlinear terms
$$ f(u_1,u_2) = (b/a) u_1^2 + 2 a u_1 u_2 + u_1^2 u_2 .$$
The long-wave Hopf instability occurs at the
critical wave number
$ k =0$  for $ b= b_{hopf}(a) = 1+ a^2 $.  For more details see \cite{SW22}.
This system can be brought into the form \eqref{rd1}-\eqref{rd2} by 
introducing a variable $ v $ satisfying 
$$ 
\partial_t v = d_v \partial_x^2  v + \partial_x^2 g(u_1,u_2) ,
$$
with $ g(0,0) = 0 $, by replacing $ b $ by $ v + \widetilde{b} $, 
and by introducing the small bifurcation parameter
$ \epsilon^2 = (\widetilde{b}-b_{hopf})/b_{hopf} $.
}
 \end{remark}
 
 The {\bf plan of the paper} is as follows. In Section \ref{secGL}
 we discuss the global existence and uniqueness of solutions of the 
 amplitude system \eqref{as1intro}-\eqref{as2intro}. The proof will be given in 
 Appendix \ref{appD}. Section \ref{secnormal} contains a number 
 of preparations, in particular we eliminate a number of oscillatory 
 terms from \eqref{rd1}-\eqref{rd2} by so called normal form transformations.
 In Section \ref{secGL4} we derive the amplitude equations and define  
 the Ginzburg-Landau manifold, the set of solutions which can be approximated 
 by our amplitude system. In Section \ref{secglobal} we formulate the attractivity result which is proven in Appendix \ref{appB}, the approximation result 
which is proven in Appendix \ref{appC} and put them together to conclude on the 
global existence and uniqueness of solutions of the original reaction-diffusion 
system \eqref{rd1}-\eqref{rd2}. In Section \ref{secdisc} a few further 
questions are discussed.
Moreover, in  Appendix \ref{appA2} some estimates are provided which are used in the sequel.
\medskip

{\bf Notation.} \label{notat}
The Sobolev space $ H^s $ is equipped with the norm
$
\| u \|_{H^s} = \sum_{j= 0}^s \| \partial_x^j u \|_{L^2}
$, where $ \| u \|_{L^2}^2 = \int |u(x)|^2 dx $.
The  space $ H^s_{l,u} $ of $ s $-times locally  uniformly weakly 
differentiable functions is equipped with the norm
$
\| u \|_{H^s_{l,u} } = \sum_{j= 0}^s \| \partial_x^j u \|_{L^2_{l,u}}
$, where $ \| u \|_{L^2_{l,u}} = \sup_{x \in \R}  (\int_x^{x+1} |u(y)|^2 dy )^{1/2}$,
cf \cite[\S 8.3.1]{SU2017book}.
Fourier transform w.r.t. the spatial variable is denoted by $ \mathcal{F} $ and the inverse 
Fourier transform  by $ \mathcal{F}^{-1} $.
Possibly different constants which can be chosen 
independently of the small perturbation parameter $ 0 < \varepsilon \ll 1 $ 
are often denoted with the same symbol $ C $.
\medskip

{\bf Acknowledgement.} The paper is partially supported by the Deutsche Forschungsgemeinschaft 
DFG under the grant Schn520/10.

\section{Analysis of the amplitude system}
\label{secGL}

We consider 
\begin{eqnarray} \label{as3}
\partial_T A  & = & a_0 \partial_X^2 A + a_1 A+ a_2   A B-  a_3  A |A|^2,\\
\partial_T B & =  &  b_0 \partial_X^2 B +  b_1 \partial_X^2 (|A|^2), \label{as4}
\end{eqnarray}
where $ T \geq 0 $, $ X \in \R $, $ A(X,T) \in \C $, $ B(X,T) \in \R $, and with
coefficients having properties as specified below the equations \eqref{as1intro}-\eqref{as2intro}.
We are interested in the situation  of an unstable trivial solution, i.e., $ a_1 > 0 $. 
This is the general form of the amplitude system which appears 
for a long wave Hopf bifurcation in a pattern forming system with a conservation law.
The system has been derived for pattern forming systems with a conservation law exhibiting 
a Turing instability, too, cf. \cite{MC00}. In a singular limit spike solutions have been constructed
in \cite{NWW02}.

\begin{remark}{\rm \label{rems}
By rescaling $ A $, $ B  $, $ T $, and $ X $ and by possibly changing the sign of $ B $,
four of the coefficients 
can be eliminated. 
We set 
$$ A= c_A   \widetilde{A} , \qquad B= c_B   \widetilde{B} ,
\qquad T= c_T   \widetilde{T}  , \qquad \textrm{and}
\qquad X= c_X   \widetilde{X}.
$$ 
We find 
\begin{eqnarray*}
\partial_{\widetilde{T}}   \widetilde{A}   & = & c_T a_0 c_X^{-2}\partial_{\widetilde{X}}^2  \widetilde{A}  + c_T a_1  \widetilde{A} + c_Ta_2  c_B  \widetilde{A}   \widetilde{B} -  c_T a_3 c_A^2   \widetilde{A}  | \widetilde{A} |^2,\\
\partial_{\widetilde{T}}   \widetilde{B}  & =  &  c_T b_0 c_X^{-2} \partial_{\widetilde{X}}^2  \widetilde{B}  +  c_T b_1  c_A^2 c_X^{-2}c_B^{-1}\partial_X^2 (| \widetilde{A} |^2).
\end{eqnarray*}
We first choose $ c_T  \in \mathbb{R}$ such that $ {c_T} a_1  = 1  $.  Next we set 
$ c_A > 0 $ such that $ c_T (\textrm{Re}a_3) c_A^2  = 1 $. Then we 
choose  $ c_X > 0 $ such that 
$ c_T (\textrm{Re}a_0) c_X^{-2} = 1 $. Finally, we set $ c_B \in \mathbb{R}$ such that 
$ c_T b_1  c_A^2 c_X^{-2}c_B^{-1} = 1 $
if $ b_1 \neq  0 $. If $ b_1 =   0 $, subsequently in \eqref{as6} the term 
$  \partial_X^2 (|A|^2) $ will be away.
Defining 
$$
\beta = c_Ta_2  c_B , \quad 
\alpha = c_T b_0 c_X^{-2} , 
\quad 
\gamma_0 = \textrm{Im}(c_T a_0 c_X^{-2}),
\quad 
\gamma_3 = \textrm{Im}( c_T a_3 c_A^2  )
 $$
and dropping the tildes we finally consider
\begin{eqnarray} \label{as5}
\partial_T A  & = & (1+i \gamma_0) \partial_X^2 A +  A+   \beta  A B- (1+i \gamma_3)   A |A|^2,\\
\partial_T B & =  &  \alpha \partial_X^2 B +  \partial_X^2 (|A|^2), \label{as6}
\end{eqnarray}
with $ \alpha> 0 $ and $ \beta, \gamma_0,\gamma_3 \in \R $.}
\end{remark}

\begin{remark}{\rm
Before we discuss the local and global existence of this system we have a short look 
at a family of special solutions.
There are the $ X $-independent time-periodic solutions 
$ B = b $, $ A = \widehat{A} e^{i \omega T}$ with $ |\widehat{A}|^2 = 1+ \beta b $ 
and $ \omega =  - |\widehat{A}|^2 \gamma_3 $ 
for every $ b $ with  $ 1+ \beta b > 0 $.
In case $ 1+ \beta b \leq  0 $ we have the stationary  solutions $ B = b $ and $ A = 0 $.
}
\end{remark}

\begin{remark} \label{rem23}{\rm
Global existence for the classical Ginzburg-Landau equation on the real line,  \eqref{as5} in case $ \beta =0 $,
can be obtained in $ C^0_b(\mathbb{R}) $  with
the maximum principle if $ \gamma_0 = \gamma_3 = 0 $. By the smoothing of the diffusion semigroup, 
global existence  follows in all $ C^n_b $-spaces and $ H^m_{l,u} $-spaces for $ m > 1/2 $.
An approach for general $ \gamma_0 $ and $ \gamma_3 $ is to work with weighted energies 
$ \int_{\mathbb{R}}  \rho_{\delta}(X) |A(X)|^2 dX $,
where  $  \rho_{\delta}(X)= (1+(\delta X)^2)^{-1} $ for $ \delta > 0 $, cf. 
\cite{MS95}.
}
\end{remark}
\begin{remark} \label{rem24}{\rm
However, so far, both approaches described in Remark \ref{rem23} do not give global existence 
for the amplitude  system  \eqref{as5}-\eqref{as6} on the real line. 
Weighted energy estimates gives via the linear terms $ \partial_X^2 A $ and 
$ \alpha \partial_X^2 B $
some exponential growth of order $ \mathcal{O}(\delta^2) $.
For the classical Ginzburg-Landau equation one can get rid of these growth rates 
with  the $ - |A|^2 A $-term which allows for a point-wise 
estimate
\begin{equation} \label{FRA}
\int_{\mathbb{R}}  \rho_{\delta}(X)(|A(X)|^2 - |A(X)|^4) dX \leq 
\int_{\mathbb{R}}  \rho_{\delta}(X)(1 - |A(X)|^2) dX .
\end{equation}
However, there is no counterpart in  \eqref{as5}-\eqref{as6} 
which can stop the growth 
of  the weighted $ B $-variable.
}
\end{remark}

We help ourselves by considering the amplitude system \eqref{as5}-\eqref{as6} 
with periodic boundary conditions.
$ 2 \pi $-periodicity for \eqref{as3}-\eqref{as4} corresponds to 
$ L $-periodicity for \eqref{as5}-\eqref{as6} with $ L =  2\pi \sqrt{a_1/a_0}$. 

\begin{remark} \label{rem20}{\rm
In case of periodic boundary conditions we could always assume that the mean value $ b $ of $ B $ vanishes. If this would not be the case, we could set $ B = b + \widetilde{B} $, with $ \widetilde{B} $ having 
a vanishing mean value. Then we would  obtain
\begin{eqnarray*} 
\partial_T A  & = & (1+ i \gamma_0) \partial_X^2 A +  A+ \beta   A (b + \widetilde{B})- (1+ i \gamma_3)   A |A|^2,\\
\partial_T \widetilde{B}  & =  &  \alpha \partial_X^2 \widetilde{B}  +  \partial_X^2 (|A|^2).
\end{eqnarray*}
Hence, by redefining the coefficient $ a_1 $ we could always come to a system, for which the
mean value of $ B $ vanishes for all $ T \geq 0 $.
}
\end{remark}

The choice of  periodic boundary conditions allows us to use classical energy estimates without weights. In case $ 1+  \alpha^{-1} \beta > 0 $ we have the 
following global existence result.

\begin{theorem}  \label{glglob}
Assume that $ 1+  \alpha^{-1} \beta > 0 $ holds.
Fix $ s \geq 0 $, $ L > 0 $ and consider \eqref{as5}-\eqref{as6} with $ L$-periodic boundary conditions.
Then there exists a  $ C_2 > 0 $ such that for all $ C_1 > 0 $ there exists a $ T_0 > 0 $ such that the following holds. For initial conditions 
$ (A(\cdot,0),B(\cdot,0)) \in H^{s+1} \times H^s $ with 
$ \int_0^L B(X,0) dX = 0 $ and  
$$ 
\|A(\cdot,0) \|_{H^{s+1} } + \|B(\cdot,0)\|_{H^{s} }  \leq C_1
$$ 
the associated unique global solution $ (A,B) \in C([0,\infty) , H^{s+1}  \times H^s )  $ 
satisfies 
$$ 
\|A(\cdot,T) \|_{H^{s+1}} + \|B(\cdot,T)  \|_{H^{s}}  \leq C_2
$$ 
for all $ T \geq T_0 $.
\end{theorem}
\noindent
{\bf Proof.} See Appendix \ref{appD}. \qed

\section{Some preparations}
\label{secnormal}

All operators appearing in the following are 
so called multipliers. 
A linear  operator $M$ is called
multiplier if there exists a function $ \widehat{M}:\R\to \C $ such that 
$ M u = \mathcal{F}^{-1} (\widehat{M} \mathcal{F} u ) $, i.e., if the associated operator  is 
a multiplication operator in Fourier space. 
Typical examples are differential operators, semigroups, or mode-filters,
but also the normal form transformations at the end of this section 
can be interpreted as  multilinear multipliers.

\subsection{The mode-filters}

For estimating the different parts of the solutions we use 
so called mode-filters. 
Since we work in $ H^n_{l,u} $-spaces
we cannot use cut-off functions in Fourier space to extract certain modes 
from the solutions. 
The associated operators in $ H^n_{l,u} $ would not be smooth and so 
we  take a $ \widehat{\chi} \in C_0^{\infty} $
with
\begin{equation} \label{chi}
\widehat{\chi}(k) = \left\{ \begin{array}{cc} 1 ,& \textrm{for } |k| \leq 0.45 \widetilde{\delta} ,\\
0 ,& \textrm{for } |k| \geq 0.55 \widetilde{\delta} ,\\
\in [0,1], &  \textrm{else},
\end{array} \right.
\end{equation}
for a $ \widetilde{\delta} > 0 $ sufficiently small but independent of the small perturbation parameter $ 0 < \varepsilon^2 \ll 1 $.
For extracting the modes around the Fourier wave number $ k= 0 $ we define 
a mode-filter $ E_0 $ by 
$$ 
\widehat{E}_0(k)  \widehat{u}(k) = \widehat{\chi}(k) \widehat{u}(k).
$$
This operator can be estimated as follows.
\begin{lemma}  \label{lem33}
For every $ m \in \mathbb{N}_0 $ the operator $ E_0 $ is a  bounded operator from $ L^2_{l,u} $ to $ H^m_{l,u} $,
in detail, there exist constants $ C_{m} $ such that $ \| E_0 \|_{L^2_{l,u}  \to H^m_{l,u} }  \leq C_{m}  $.  
%
\end{lemma}
\noindent
{\bf Proof.} We use multiplier theory in $ H^m_{l,u} $-spaces, cf. \cite[\S 8.3.1]{SU2017book}. We have 
$$ \| E_0 v \|_{H^m_{l,u}} \leq C \| \widehat{M} \|_{C^2_b}\| v \|_{L^2_{l,u}},$$
with 
$
\widehat{M}(k) = (1+k^2)^{m/2} \widehat{\chi}(k)
$. 
\qed

\subsection{The normal form transformation}
\label{secnft}

For the subsequent analysis we need a  separation of the 
$ u $-modes into exponentially damped  $ (j=3,\ldots,d) $ and critical modes $ (j=\pm 1) $.
In order to do so, we let
$$ \tilde{P}_{\pm 1} (k,\epsilon^2){u} = \frac{1}{2 \pi i }
\int_{\Gamma_{\pm 1}} ( \mu Id. - \tilde{\Lambda}_u
(k,\epsilon^2) )^{-1} {u} d\mu ,
$$
where $  \tilde{\Lambda}_u = \mathcal{F} \Lambda_u \mathcal{F}^{-1}$, with $ \Lambda_u $ defined in \eqref{lrd1},
and where
$ \Gamma_{\pm 1} $ is a closed curve surrounding
the single eigenvalue $ \lambda_{\pm 1}|_{\varepsilon = 0,k=0}= \pm i \omega_0 $ anti-clockwise.
By the assumption (Spec) the projections $  \tilde{P}_j $ can be defined
for wave numbers in a neighborhood $ U_\rho(0) $ for a $ \rho > 0 $
and so we set 
$$ 
E_{\pm 1} = E_0  \tilde{P}_{\pm 1} , \quad  E_c = E_1 + E_{-1}, \qquad \textrm{and} \qquad E_s = Id.-E_{c},
$$
choosing $\widetilde{\delta} < \rho/2 $ in \eqref{chi}.
Moreover, we define scalar-valued projections $ \tilde{p}_{\pm 1} $ by 
$$ 
\tilde{P}_{\pm 1} (k,\epsilon^2)  u = (\tilde{p}_{\pm 1} (k,\epsilon^2) u)  \widehat{U}_{\pm 1}(k,\epsilon^2)
$$ 
and $ e_{\pm 1} = E_0  \tilde{p}_{\pm 1} $.
With these operators we separate our system \eqref{rd1}-\eqref{rd2} 
in critical, neutral,  and exponentially damped modes.

Then, in Fourier space,  we write 
$$ 
\widehat{u}(k,t) = \widehat{c}_1(k,t) \widehat{U}_1(k) + \widehat{c}_{-1}(k,t) \widehat{U}_{-1}(k)
+ \widehat{u}_s(k,t),
$$ 
with $ \widehat{c}_{\pm 1}(k,t) \in \C $, and define $ c_{\pm 1} $ and $ u_s $ 
to be solutions of 
\begin{eqnarray} \label{rdnf1}
\partial_t c_1 & = & \lambda_1 c_1 + f_{1}(c_1,u_s,v) , \\
\partial_t c_{-1} & = & \lambda_{-1} c_{-1} + f_{-1}(c_1,u_s,v) , \\
\partial_t u_s & = & \Lambda_s u_s + f_{s}(c_1,u_s,v) ,\\
\partial_t v & = & \Lambda_v v  + \partial_x^2 g(c_1,u_s) ,\label{rdnf3}
\end{eqnarray}
with 
\begin{eqnarray*} 
f_{\pm 1}(c_1,u_s,v)  & = & e_1 f(u,v) = \mathcal{O}(|c_1|^2+ |u_s|^2+(|u_1|+|u_s|) |v|) , \\
f_{s}(c_1,u_s,v) & = & E_s f(u,v) =  \mathcal{O}(|c_1|^2+ |u_s|^2+(|u_1|+|u_s|) |v|),\\
g(c_1,u_s)  & = & \mathcal{O}(|c_1|^2+ |u_s|^2).
\end{eqnarray*}
Since $ c_{-1} = \overline{c_1} $ we do not explicitly denote the appearance of 
$ c_{-1} $  in
various places.

Since $ c_{1} $ approximately oscillates as  $ e^{ i \omega_0 t} $ all quadratic combinations of $ c_{1} $ and $ c_{-1} $ can be eliminated 
from the $ c_{1} $-equations by a near identity change of variables 
$$
\check{u}_1 = c_1 + \mathcal{O}(|c_1|^2).
$$ 
A similar statement holds for the $ c_{-1} $-equation.
For details see the subsequent Remark \ref{rem32sec3}.
\begin{remark}{\rm
In a similar way  terms $ v  c_{\pm 1}$ in the $ v $-equation could be eliminated in case of a more general nonlinearity in the $ v $-equation.}
\end{remark}

After the transform we have a system of the form 
 \begin{eqnarray} \label{imag1}
\partial_t \check{u}_1 & = & \lambda_1 \check{u}_1 + 
\check{f}_1(\check{u}_1,u_s,v)
 , \\
\partial_t u_s & = & \Lambda_s u_s +\check{f}_s(\check{u}_1,u_s,v),\\
\partial_t v & = & \Lambda_v v  + \partial_x^2 \check{g}(\check{u}_1,u_s), \label{imag3}
\end{eqnarray}
with 
\begin{eqnarray*} 
\check{f}_1(\check{u}_1,u_s,v) & = & \mathcal{O}(|\check{u}_1|^3+ |\check{u}_1||u_s| +  |u_s|^2+(|\check{u}_1|+|u_s|) |v|) , \\
\check{f}_s(\check{u}_1,u_s,v) & = & \mathcal{O}(|\check{u}_1|^2+ |u_s|^2+(|\check{u}_1|+|u_s|) |v|),\\
\check{g}(\check{u}_1,u_s) & = & \mathcal{O}(|\check{u}_1|^2+ |u_s|^2).
\end{eqnarray*}

Detailed estimates about this transformation and the nonlinear terms are given below when needed.
\begin{remark}
\label{rem32sec3}
{\rm
In lowest order the equation for $ c_{1} $ is of the form
$$ 
\partial_t c_{1}  = \lambda_1 c_{1}   + N_{1,1} (c_{1} ,c_{1} )
+ N_{1,-1} (c_{1} ,c_{-1} )+ N_{-1,-1} (c_{-1} ,c_{-1} ) + h.o.t.
$$ 
where in Fourier space the $ N_{i,j} $ have a representation 
$$ \hat{N}_{i,j}(c_{i} ,c_{j} )[k] = \int \hat{n}_{i,j}(k,k-m,m) \hat{c}_{i} (k-m) \hat{c}_{j}(m) dm ,$$
with kernel functions $ \hat{n}_{i,j} : \R^3 \to \C $.
The quadratic terms can be eliminated by a transform 
$$
\check{u}_{1}  =  c_{1}   + B_{1,1} (c_{1} ,c_{1} )
+ B_{1,-1} (c_{1} ,c_{-1} )+ B_{-1,-1} (c_{-1} ,c_{-1} )
$$
where in Fourier space the $ B_{i,j} $ have a representation 
$$ \hat{B}_{i,j}(c_{i} ,c_{j} )[k] = \int \hat{b}_{i,j}(k,k-m,m) \hat{c}_{i} (k-m) \hat{c}_{j}(m) dm .$$
The kernels $ \hat{b}_{i,j}(k,k-m,m) $ are solutions of 
$$ 
(\tilde{\lambda}_{1}(k) - \tilde{\lambda}_{i}(k-m) -
\tilde{\lambda}_{j}(m))\hat{b}_{i,j}(k,k-m,m) = \hat{n}_{i,j}(k,k-m,m) 
$$
which are well-defined and bounded since 
$$
\inf_{k,m \in U_{4 \rho}(0)} |
\widehat{\lambda}_{j_1}(k) - \widehat{\lambda}_{j_2}(k-m) -
\widehat{\lambda}_{j_3}(m)| \geq C > 0
$$
for all $ j_1,j_2,j_3 \in \{-1,1\} $.
For more details see \cite[\S 11]{SU2017book} or  \cite[\S 4]{Schn98JNS}.
}\end{remark}
\begin{remark}{\rm
After the transform we have a system of the form 
\begin{eqnarray*}
\partial_t \check{u}_{1} & = & \lambda_1 \check{u}_{1}   + N_{1,1,1} (\check{u}_{1} ,\check{u}_{1} ,\check{u}_{1})
+ N_{1,1,-1} (\check{u}_{1} ,\check{u}_{1},\check{u}_{-1} ) \\ && + N_{1,-1,-1} (\check{u}_{1},\check{u}_{-1} ,\check{u}_{-1} ) +
N_{-1,-1,-1} (\check{u}_{-1},\check{u}_{-1} ,\check{u}_{-1} ) 
+ h.o.t.
\end{eqnarray*}
where in Fourier space the $ N_{i,j,k} $ have a similar 
representation as above.
Except of  $ N_{1,1,-1} $
the three other terms are non-resonant such that these can be eliminated 
by a second transformation. 
}\end{remark}
\begin{example}{\rm 
Applying the normal form transformation to the system from Example 
\ref{exam11} yields a system of the form 
\begin{eqnarray*}
\partial_t u_1 & = & \partial_x^2 u_1 + i \omega_0 u_1 + \varepsilon^2 u_1 
  + v u_1  - \gamma_3 u_1^2 u_{-1}  + h.o.t.\\
\partial_t u_{-1}  & = & \partial_x^2 u_{-1}  - i \omega_0 u_{-1}  + \varepsilon^2 u_{-1}  
 + v u_{-1}-  \overline{\gamma_3}  u_{-1}^2 u_{1}+ h.o.t., \\
\partial_t v & = & \partial_x^2 v +   \partial_x^2 ( u_1 u_{-1})+ h.o.t. ,
\end{eqnarray*}
with $ \gamma_3 $ given by \eqref{gamma3}.}
\end{example}

\section{The Ginzburg-Landau manifold}
\label{secGL4}

The notation Ginzburg-Landau manifold or Ginzburg-Landau set, cf. \cite{Eck93},
was chosen to describe the set of initial conditions of the original system 
\eqref{rd1}-\eqref{rd2} which can be described 
by the Ginzburg-Landau approximation.
In the non-conservation law case it was shown that this set is attractive,
cf. \cite{Eck93,BSvH94,Schn95c}. In the  conservation law case
a first result was established in \cite{DSSZ16}. We will come back to this 
in Section \ref{secatt}. It is the purpose of this section to derive
the amplitude system, to compute a higher order approximation and to define 
what we will mean by Ginzburg-Landau manifold.

For possible future applications, similar to \cite{MS95,Schn99JMPA}, we introduce a new perturbation parameter $ \delta  $ with $ 0 < \varepsilon \leq \delta \ll 1 $ and distinguish this parameter 
from the bifurcation parameter $ 0 < \varepsilon  \ll 1 $. 

\subsection{Derivation of the amplitude system}

Our starting point for the derivation of the amplitude system is System \eqref{rdnf1}-\eqref{rdnf3} which we write as 
\begin{eqnarray*} 
\textrm{Res}_1 & = &-\partial_t c_1 + \lambda_1 c_1 + f_1(c_1,u_s,v), \\
\textrm{Res}_s & = &-\partial_t u_s  +\Lambda_s u_s + f_s(c_1,u_s,v),\\
\textrm{Res}_v & = &-\partial_t v  +\Lambda_v v  + \partial_x^2 g(c_1,u_s).
\end{eqnarray*}
The so called residuals $\textrm{Res}_1  $, $\textrm{Res}_s  $, 
and $\textrm{Res}_v  $ contain all terms which remain after inserting 
an approximation into  System \eqref{rdnf1}-\eqref{rdnf3}.

For the derivation of the amplitude system, cf. Example \ref{exam11},
we need an  ansatz 
\begin{eqnarray*}
c_{1}(x,t) & = & \delta A_1(X,T) e^{i \omega_0 t}
+ \delta^2 A_{1,0}(X,T) \\
&&
+ \delta^2 A_{1,2}(X,T) e^{2 i \omega_0 t} + \delta^2 A_{1,-2}(X,T) e^{-2 i \omega_0 t}
, \\
c_{-1}(x,t) & = & \delta A_{-1}(X,T) e^{-i \omega_0 t}
+ \delta^2 A_{-1,0}(X,T) \\
&& 
+ \delta^2 A_{-1,2}(X,T) e^{2 i \omega_0 t}
+ \delta^2 A_{-1,-2}(X,T) e^{-2 i \omega_0 t} 
, \\
u_s (x,t) & = &  
\delta^2 A_{s,2}(X,T) e^{2 i \omega_0 t}
+ \delta^2 A_{s,0}(X,T) 
+ \delta^2 A_{s,-2}(X,T) e^{-2 i \omega_0 t} ,\\
v(x,t) & = & \delta^2 B_0(X,T),
\end{eqnarray*}
with $ X = \delta x $ and $ T = \delta^2 t $.
By equating the coefficients in front of $ \delta^2 e^{in\omega_0 t} $, with $ n = 0 , \pm 2 $, to zero, 
we find $ A_{j,2},  A_{j,0}, A_{j,-2} $ for $ j = -1,1,s $  
as solutions of equations of the form 
\begin{eqnarray*}
A_{j,2} & = & \gamma_{j,2} A_{1} A_{1} , \\  
A_{j,0} & = & \gamma_{j,0} A_{1} A_{-1} , \\  
A_{j,-2} & = & \gamma_{j,-2} A_{-1} A_{-1} ,
\end{eqnarray*}
with coefficients  $ \gamma_{j,i} $.
The $ A_1 $, $ A_{-1} $, and $ B_0 $ satisfy a system of the form 
\begin{eqnarray*}
\partial_T A_1  & = & a_0 \partial_X^2 A_1 + \frac{\varepsilon^2}{\delta^2}a_1 A_1+ a_2   A_1 B_0-  \widetilde{a}_4  A_1 |A_1|^2 \\ && + \sum_{j = \pm 1, s} a_{6,j} A_{j,2} A_{-1} + \sum_{j = \pm 1, s} a_{7,j}  A_{j,0} A_{1} ,\\
\partial_T B_0 & =  &  b_0 \partial_X^2 B_0 +  b_1 \partial_X^2 (|A_1|^2), 
\end{eqnarray*}
Eliminating the $ A_{j,2},  A_{j,0}, A_{j,-2} $ for $ j = -1,1,s $  through the above equations gives the amplitude system 
\begin{eqnarray} \label{as1introg}
\partial_T A_1  & = & a_0 \partial_X^2 A_1 + \frac{\varepsilon^2}{\delta^2} a_1 A_1+  a_2  A_1 B_0 -  a_3  A_1 |A_1|^2,\\
\partial_T B_0& =  &  b_0 \partial_X^2 B_0 +  b_1 \partial_X^2 (|A_1|^2), \label{as2introg}
\end{eqnarray}
similar to 
\eqref{as1intro}-\eqref{as2intro}.
 We formally have 
 $$\textrm{Res}_1 = \mathcal{O}(\delta^3) , \quad \textrm{Res}_s= \mathcal{O}(\delta^3)  , \quad  
 \textrm{Res}_v = \mathcal{O}(\delta^4)  $$
for this approximation. In the residual of the $ c_1 $-equation we have for instance a term $ \delta^3 A_1^3 e^{3 i \omega_0 t} $ and  in the residual of the $ v $-equation we have for instance a term $ \delta^4 \partial_X^2 (A_1^2)  e^{2 i \omega_0 t} $.

In order to show that the amplitude system \eqref{as1introg}-\eqref{as2introg}
makes correct predictions about the original system 
\eqref{rd1}-\eqref{rd2} we establish subsequently the approximation
Theorem \ref{thapp}.

\subsection{Construction of a higher order approximation}
\label{secgl42}

In order to obtain a more precise approximation we add higher order terms 
to the previous approximation. We insert 
$$
c_{1} = \psi_{1}, \qquad c_{-1} = \psi_{-1}, \qquad u_{s} = \psi_{s}
\qquad v = \psi_{v}
$$
with 
\begin{eqnarray*}
\psi_{1}(x,t) & = & \sum_{m = -N}^N \sum_{n=0}^{M_1(N,m)} \delta^{\beta_{1}(m)+n} A_{+,m,n}(X,T) e^{ i m \omega_0 t}
, \\
 \psi_{-1}(x,t) & = &  \sum_{m = -N}^N \sum_{n=0}^{M_1(N,m)} \delta^{\beta_{-1}(m)+n} A_{-,m,n}(X,T) e^{ i m \omega_0 t}
, \\
\psi_{s} (x,t) & = &  
\sum_{m = -N}^N \sum_{n=0}^{M_s(N,m)} \delta^{\beta_{s}(m)+n} A_{s,m,n}(X,T) e^{ i m \omega_0 t},\\
\psi_{v}(x,t) & = & \sum_{m = -N}^N \sum_{n=0}^{M_v(N,m)} \delta^{\beta_{v}(m)+n} B_{m,n}(X,T) e^{ i m \omega_0 t},
\end{eqnarray*}
where $ N $, $ M_1(N,m) $, $ M_s(N,m) $, and $ M_v(N,m) $ are sufficiently large 
numbers 
such that 
$$\textrm{Res}_c = \mathcal{O}(\delta^{\theta+2}) , \quad \textrm{Res}_s= \mathcal{O}(\delta^{\theta+2})  , \quad  
 \textrm{Res}_v = \mathcal{O}(\delta^{\theta+2})  $$
for a given $ \theta \in \mathbb{N}$
and where 
$$
\begin{array}{|c||c|c|c|c|c|c|c|c|}\hline
 m  &  - 3 & - 2 & -1 & 0 & 1 & 2 & 3 & m \\ \hline \hline
\beta_{1}(m)  & 3 & 2 & 3 & 2 & 1 & 2 & 3 & m \\ \hline
\beta_{-1}(m) & 3 & 2 & 1 & 2 & 3 & 2 & 3 & m\\ \hline
\beta_{s}(m)  & 3 & 2 & 3 & 2 & 3 & 2 & 3 & m\\ \hline
\beta_{v}(m)  & 5 & 4 & 5 & 2 & 5 & 4 & 5 & m+2 \\ \hline 
\end{array}
$$
The associated approximation is then denoted with $ \Psi_\theta $.

The coefficient functions are determined as follows.
The functions $ A_{+,1,0} $, $ A_{-,-1,0} $, and $ B_{0,0} $ satisfy the amplitude system from above.
 The $ A_{+,1,n} $, $ A_{-,-1,n} $, and $ B_{0,n} $ for $ n \geq 1 $ 
 satisfy linearisations of the amplitude system from above with some inhomogeneous 
 terms which in the end depend on terms $ A_{+,1,j} $, $ A_{-,-1,j} $, and $ B_{0,j} $ for $ 0 \leq j \leq n-1 $.  All other $ A_{+,m,n} $, $ A_{-,m,n} $, $ A_{s,m,n} $, and $ B_{m,n} $
 satisfy algebraic equations and can be computed in terms of the 
 $ A_{+,1,j} $, $ A_{-,-1,j} $, and $ B_{0,j} $ for $ 0 \leq j \leq n $. 
 
 The solutions of this system are uniquely determined by the set of initial conditions 
 $ A_{+,1,j}|_{T=0} $, $ A_{-,-1,j}|_{T=0} $, and $ B_{0,j}|_{T=0} $ for $ 0 \leq j \leq n $.
 
 \begin{definition}
For initial conditions 
$$
A_{+,1,0}|_{T=0}  = A_1|_{T=0}  ,  \qquad  A_{-,-1,0}|_{T=0}  = \overline{A_1}|_{T=0}  ,
\qquad  B_{0,0}|_{T=0} = B_0|_{T=0}  
$$
and 
$$ A_{+,1,j}|_{T=0} ,  \qquad  A_{-,-1,j}|_{T=0} , \qquad B_{0,j}|_{T=0} $$ 
determined by the construction in Appendix \ref{appb2}
for $ 1 \leq j \leq n $ and $ (A_1,B_0) $ satisfying \eqref{as1introg}-\eqref{as2introg} we call the set of approximate solutions 
$$ 
(u,v)(\cdot,t) = \Psi_\theta(A_1(\cdot,T) ,B_0(\cdot,T))
$$ 
for the original system \eqref{rd1}-\eqref{rd2}
the 
Ginzburg-Landau manifold, 
where $ \Psi_\theta $ is  the associated higher order approximation defined above. 
\end{definition}

\section{The global existence and uniqueness result}
\label{secglobal}

Throughout the rest of this paper 
we replace the boundary conditions \eqref{uvper} by the  boundary conditions 
\begin{equation} \label{uvdeltaper}
u(x,t) = u(x+2 \pi/\delta,t) \qquad \textrm{and} \qquad 
v(x,t) = v(x+2 \pi/\delta,t) .
\end{equation} 
with $ 0 < \varepsilon \leq \delta \ll 1 $ and set later on $ \delta=  \varepsilon$.
\begin{remark}\label{rem51}{\rm
There is  local existence and uniqueness of (mild) solutions 
$$ 
(u,v) \in C([0,t_0],H^{n+1}_{l,u} \times  H^{n}_{l,u} )
$$ 
of 
\eqref{rd1}-\eqref{rd2}  for initial conditions $ (u_0,v_0) \in H^{n+1}_{l,u} \times  H^{n}_{l,u}  $ 
if $ n \geq  0 $
where the existence time $ t_0 > 0 $ only depends on 
$ \| u_0 \|_{H^{n+1}_{l,u} } + \| v_0 \|_{H^{n}_{l,u} }  $. 
This can be established with the standard  fixed point argument
for semilinear  parabolic equations, cf. \cite{He81}. For $ n \geq 0 $  the right-hand side
 of the variation of constant formula associated to \eqref{rd1}-\eqref{rd2}  is a contraction in a ball 
 in $ C([0,t_0], H^{n+1}_{l,u} \times  H^{n}_{l,u}) $ for $ t_0 > 0 $ sufficiently small 
 using that the nonlinear terms $ (f(u,v),\partial_x g(u))  $ are 
smooth mappings from $  H^{n+1}_{l,u} \times  H^{n}_{l,u} $ 
to $ H^{n}_{l,u} \times  H^{n}_{l,u} $ and that the linear semigroups 
$ (e^{D \partial_x^2 t}, e^{d_v \partial_x^2 t} \partial_x) $  map
$ H^{n}_{l,u} \times  H^{n}_{l,u} $ to $  H^{n+1}_{l,u} \times  H^{n}_{l,u} $ with an integrable  singularity
$ t^{-1/2} $. 
}\end{remark}
Hence, for 
establishing 
the global existence and uniqueness of (mild) solutions we need 
to bound the solutions in $ H^{n+1}_{l,u} \times  H^{n}_{l,u} $, i.e.,
if we establish an a priori  bound 
\begin{equation} \label{apriori}
\sup_{t \in [0,\infty)} (\| u(t) \|_{H^{n+1}_{l,u} } + \| v(t) \|_{H^{n}_{l,u} }) \leq C_3 < \infty,
\end{equation}
where $ C_3 $ is a constant only depending on $ \| u_0 \|_{H^{n+1}_{l,u} } + \| v_0 \|_{H^{n}_{l,u} }  $, then
 the 
 local existence and uniqueness theorem can be applied again and again 
and  the local solutions can be continued to global solutions
 $$ 
(u,v) \in C([0,\infty),H^{n+1}_{l,u}  \times  H^{n}_{l,u} ).
$$ 
The necessary a-priori estimates \eqref{apriori} for \eqref{rd1}-\eqref{rd2}  
can be obtained in a sufficiently small $ \mathcal{O}(\delta)$-neighborhood of the weakly unstable origin with the help of an attractivity  and  approximation result for  the Ginzburg-Landau manifold  and the existence of an absorbing ball for the amplitude 
system.

The attractivity theorem is as follows
\begin{theorem} \label{thatt}
For all $ R_0 > 0  $, $ n \geq 0  $, and all $ \theta \in \mathbb{N}_0 $ the following holds. Consider \eqref{rd1}-\eqref{rd2} 
with initial conditions $ (u_0,v_0) \in H^{n+1}_{l,u}  \times  H^{n}_{l,u}$ 
satisfying
$$ 
\| u_0 \|_{H^{n+1}_{l,u} } +\delta^{-1} \| v_0 \|_{H^{n}_{l,u}} \leq R_0 \delta.
$$ 
Then there exists 
a time $ T_1 \in (0,1) $, a $ \delta_1 > 0 $, an $ R_1 > 0 $  and a $ C_1 > 0 $, 
all only depending on $ R_0 $, $ \theta $, and $ n $,
such that 
for all $ \delta \in (0, \delta_1) $, all $ \varepsilon \in (0,\delta] $, and all $ m>1/2 $
there are $ (A_1(\cdot,0),B_0(\cdot,0)) \in H^{m+1}_{l,u}  \times  H^{m}_{l,u}$ with 
$$ 
\| A_1(\cdot,0) \|_{H^{m+1}_{l,u} } + \| B_0(\cdot,0) \|_{H^{m}_{l,u}} \leq R_1 
$$ 
such that the solution $ (u,v) $, with the 
initial conditions $ (u_0,v_0) $, satisfies at a time $ t = T_1/\delta^2 $ that 
$$ 
\|(u,\delta^{-1} v)|_{t = T_1/\delta^2} - ( \Psi_{\theta,u}, \Psi_{\theta,v})(A_1(\cdot,0),B_0(\cdot,0))
\|_{H^{m+1}_{l,u}  \times  H^{m}_{l,u}} \leq C \delta^\theta.
$$ 
\end{theorem}
\noindent {\bf Proof.} 
See Appendix \ref{appB}. \qed
\medskip

The dynamics on the 
Ginzburg-Landau manifold is  determined by the amplitude system 
\eqref{as1introg}-\eqref{as2introg}.
Although the 
Ginzburg-Landau manifold, 
constructed above, is not invariant under the flow of the original 
system \eqref{rd1}-\eqref{rd2}, it  is a good approximation  of the flow near the 
Ginzburg-Landau manifold.
This is documented in the following approximation theorem.

\begin{theorem} \label{thapp}
For all $R_2,T_0,C_2 > 0 $, $ n \geq 0 $  and all $ \theta \in \mathbb{N}_0 $ 
there exists $ C_3,\delta_0 >0 $ and $ m \geq 0 $ such that for all $ 0 \leq
\varepsilon \leq\delta\leq
\delta_0 $ the following holds: Let
$ (A_1,B_0)  $
be a solution of \eqref{as1introg}-\eqref{as2introg}  with 
$$ \sup_{T \in [0,T_0]}(
\| A_1(\cdot,T) \|_{H^{m+1}_{l,u}  } + \| B_0(\cdot,T) \|_{H^{m}_{l,u} } )\leq R_2 ,$$ 
with initial conditions $ (A_1,B_0)|_{T=0} = (A_1(\cdot,0),B_0(\cdot,0)) $,
and $ (u_0,v_0) \in H^{n+1}_{l,u}   \times H^{n}_{l,u}  $ with
$$ 
\|(u_0, \delta^{-1} v_0) - ( \Psi_{\theta,u}, \Psi_{\theta,v})(A_1(\cdot,0),B_0(\cdot,0))
\|_{H^{n+1}_{l,u}   \times H^{n}_{l,u} } \leq C_2 \delta^\theta.
$$  
Then
there exists a solution $ (u,v) $ of \eqref{rd1}-\eqref{rd2}  with
initial condition
$ (u,v)|_{t=0} = (u_0,v_0)  $ and
$$
\sup_{0 \leq t \leq T_0/\delta^2} \| (u,\delta^{-1} v)(\cdot,t) - ( \Psi_{\theta,u}, \Psi_{\theta,v})(A_1,B_0)(\cdot,t)
\|_{H^{n+1}_{l,u}   \times H^{n}_{l,u} } \leq C_3 \delta^\theta.
$$
\end{theorem}
\noindent {\bf Proof.} 
See Appendix \ref{appC}. \qed
\medskip

Now we have all ingredients for establishing a global existence result 
through some a-priori  bound \eqref{apriori}.
For $ \theta \geq 3 $ the following holds:

a) We start with the {\bf attractivity},  cf.  Theorem \ref{thatt}. 
For a sufficiently large $ R_0 > 0 $ 
we obtain 
$ R_1 > 0 $,  $ T_1 > 0 $, and $ A_1(\cdot,0) $ and $ B_0(\cdot,0) $ with 
 $$ 
\| A_1(\cdot,0) \|_{H^{m+1}_{l,u}  } + \| B_0(\cdot,0) \|_{H^{m}_{l,u}  } \leq R_1 
$$ 
such that the solution $ (u,v) $, with the 
initial conditions $ (u_0,v_0) $,  satisfies
$$ 
\|(u,\delta^{-1} v)|_{t = T_1/\delta^2} - (\Psi_{\theta,u}, \Psi_{\theta,v})(A_1(\cdot,0),B_0(\cdot,0))
\|_{H^{n+1}_{l,u}   \times H^{n}_{l,u} } \leq C \delta^\theta
$$ 
for $ \delta > 0 $ sufficiently small.

b) According to Theorem \ref{thGL} and Remark \ref{rem13n}, 
in case of periodic boundary 
conditions \eqref{ABper},  the amplitude system \eqref{as1introg}-\eqref{as2introg} possesses an {\bf absorbing ball} of radius $ C_R $ in $ H^{m+1}_{l,u}   \times H^{m}_{l,u}  $.
Solutions of
 \eqref{as1introg}-\eqref{as2introg} starting in the ball of the above  radius $ R_1 $ need a time $ T_0 $ to come to the 
absorbing ball of radius $ C_R $. 

c) We have to make sure that the original ball $ R_0 \delta $ 
for the original reaction-diffusion system 
\eqref{rd1}-\eqref{rd2} is so big  
that the Ginzburg-Landau embedding of the absorbing ball
for the amplitude system \eqref{as1introg}-\eqref{as2introg} of 
radius $ C_R $ is contained in this ball. In detail, 
for $A_1 $ and $ B_0 $ satisfying 
$$ 
\| A_1(\cdot,T_0) \|_{H^{m+1}_{l,u}  } + \| B_0(\cdot,T_0) \|_{H^{m}_{l,u} } \leq C_R
$$ 
we need that the starting radius $ R_0 $ is so big that  
$$ 
\|  ( \Psi_{\theta,u}, \Psi_{\theta,v})(A_1,B_0)(\cdot,T_0/\delta^2)
\|_{H^{n+1}_{l,u}   \times H^{n}_{l,u} } \leq R_0 \delta/2.
$$ 
d) Finally we use the {\bf approximation} property, i.e., that 
the amplitude systems \eqref{as1introg}-\eqref{as2introg}  makes 
correct predictions about the dynamics of the original system,
cf. Theorem \ref{thapp}.  Then the triangle inequality
guarantees that 
\begin{eqnarray*}
\lefteqn{\| (u,\delta^{-1} v)|_{(T_1+T_0)/\delta^2} \|_{H^{n+1}_{l,u}   \times H^{n}_{l,u} }} \\
& \leq & \|  (\Psi_{\theta,u}, \Psi_{\theta,v})(A_1,B_0)(\cdot,T_0/\delta^2)
\|_{H^{n+1}_{l,u}   \times H^{n}_{l,u}} \\ && + 
\sup_{0 \leq t \leq T_0/\delta^2} \| (u,\delta^{-1} v)(\cdot,T_1/\delta^2 + t) - ( \Psi_{\theta,u}, \Psi_{\theta,v})(A_1,B_0)(\cdot,t)
\|_{H^{n+1}_{l,u}   \times H^{n}_{l,u}}  \\ 
& \leq & R_0 \delta/2 + C_3 \delta^\theta \leq 3 R_0 \delta/4
\end{eqnarray*}
for $ \delta > 0 $ sufficiently small.
Thus, after a time $ (T_1+T_0)/\delta^2$ the flow of the original reaction-diffusion system \eqref{rd1}-\eqref{rd2}
has mapped the rescaled initial ball of radius $ R_0 \delta $ 
into the smaller rescaled ball of radius $ 3 R_0 \delta/4 $. Since the magnitude of the solution $ (u,v) $
is also controlled between $ t = 0  $ and $ t = (T_1+T_0)/\delta^2$ by our estimates, we established an a priori bound \eqref{apriori}.
Thus, with the above 
arguments the global existence and uniqueness of the solutions 
of \eqref{rd1}-\eqref{rd2} follows for $ \delta > 0 $ sufficiently small.

\begin{remark} \label{rem33}{\rm 
We remark that from a technical point of view, 
in contrast to previous approaches,
we moved the first step of the 
approximation result as stated \cite{Schn94c,MS95,Schn99JMPA} to 
the attractivity result. This allows us to combine the attractivity and approximation result more easily.
}
\end{remark}

\section{Discussion}

\label{secdisc}

Before we give proofs of the attractivity theorem \ref{thatt}, 
the approximation theorem \ref{thapp}, and of Theorem \ref{glglob}
we would like to close the paper by discussing two other points, namely the restriction to a nonlinearity  $ g = g(u) $ and the global existence question in case that  the periodic boundary conditions \eqref{ABper}  and \eqref{uvper} are dropped.
\begin{remark}
{\rm 
For us,  \eqref{rd1}-\eqref{rd2} is a toy model which already contains many features 
which are  relevant for the global existence question addressed in this paper.
The major restriction of our model \eqref{rd1}-\eqref{rd2} 
seems to be the assumption that $g(u) =   \mathcal{O}(|u|^2)  $
only depends on $ u $. 
However, 
an additional  dependence on $ v $ without further 
smoothing would lead to a quasilinear system and to functional analytic difficulties having to do with the quasilinearity of such a system, but 
not with the question addressed in this paper.
Alternatively, instead of \eqref{rd2},  one could consider 
the following semilinear toy problems 
$$ 
\partial_t v= d_v \partial_x^2 v + \partial_x^2 (1- \partial_x^2)^{-1} g(u,v) 
\qquad \textrm{or} \qquad 
\partial_t v= - \partial_x^4 v + d_v \partial_x^2 v + \partial_x^2  g(u,v) ,
$$ 
with $g(u,v) =   \mathcal{O}(|u|^2+ |v|^2)  $.
Since we are not interested in the sideband
unstable situation in the $ v $-equation at the wave number $ k = 0 $, cf. \cite{Eck65}, 
in these alternative models for  
 notational simplicity we would assume $g(u,v) =   \mathcal{O}(|u|^2+ |v|^2))  $ instead of  $g(u,v) =   \mathcal{O}(|u|^2+|v|)  $.
It is essential to remark that, w.r.t. the scaling used above, a term $ |v|^2 $ is 
of higher order than  a term $ |u|^2 $ and will not appear in the amplitude system 
\eqref{as1introg}-\eqref{as2introg}.
In hydrodynamical applications the quasilinearity of the problem 
often cannot be avoided, cf. \cite{Zi14}.
Global existence by the above approach is a problem which is unsolved 
in quasilinear situations even without a conservation law so far.
 }
\end{remark}

\begin{remark}
{\rm 
In this remark we would like to discuss a few observations about the 
global existence problem if the periodic boundary conditions 
 \eqref{ABper}  and \eqref{uvper}
are dropped.
We consider the situation when in lowest order in \eqref{as1intro}-\eqref{as2intro}
the $ B $-equation decouples from the $ A $-equation, i.e., $ b_1 = 0 $.
In this case the amplitude system in normal form is given by 
\begin{eqnarray} \label{as56}
\partial_T A  & = & (1+i \gamma_0) \partial_X^2 A +  A+   \beta  A B- (1+i \gamma_3)   A |A|^2,\\
\partial_T B & =  &  \alpha \partial_X^2 B . \label{as66}
\end{eqnarray}
By the maximum principle $ B $ stays bounded and for $ A $
a uniform bound in time can be established with the weighted energy method
explained in Remark \ref{rem23}.
Hence, the solutions of the amplitude system exist globally in time 
and stay uniformly bounded. However, due to the $ B $-equation the system does not possess an absorbing 
ball.

Adding the higher oder terms to the 
$ B $-equation gives a system  of the form 
$$ 
\partial_T B = d_v \partial_X^2 B +  \partial_X^2 (\mathcal{O}(\varepsilon)).
$$ 
With the variation of constant formula we obtain 
$$ 
B(T) = e^{d_v \partial_X^2 T} B(0) + \int_0^T e^{d_v \partial_X^2 (T-\tau)} \partial_X^2 (\mathcal{O}(\varepsilon)) d \tau
$$ 
and using the estimate 
$$
\| e^{d_v \partial_X^2 T} \partial_X^{2-2  \vartheta} \|_{H^n_{l,u} \to H^n_{l,u}}
\leq C T^{ \vartheta-1}
$$
we  expect 
$$ 
B(T) - e^{d_v \partial_X^2 T} B(0)  =   \mathcal{O} (\varepsilon \int_0^T (T-\tau)^{\vartheta-1} d\tau) =  \mathcal{O} (\varepsilon T^{\vartheta} ) 
$$
for a $ \vartheta > 0 $ arbitrarily small, but fixed.
Hence, in the $ A $-equation the term $ A B $ can grow as 
$ \mathcal{O}(\varepsilon T^{\vartheta} ) A $. It can be expected that it can be balanced with 
the $ - A |A|^2 $-term as long $ \varepsilon T^{\vartheta}  \leq \mathcal{O}(1) $, i.e. 
for $ T \leq  \mathcal{O}(\varepsilon^{-1/\vartheta}) $, i.e., on arbitrary 
long, but fixed, time scales w.r.t. $ \varepsilon $.  With more advanced estimates  we even would obtain 
$$ 
B(T) - e^{d_v \partial_X^2 T} B(0)  =   \mathcal{O} (C+\varepsilon \int_0^{T-1} (T-\tau)^{-1} d\tau) =  \mathcal{O} (C+\varepsilon \ln T ) ,
$$
respectively $ T \leq  \mathcal{O}(\exp(1/\varepsilon)) $. 
It will be the subject of future research to  make these arguments rigorous 
by iterating the attractivity and approximation result for a growing sequence of perturbation parameters $ \delta $. Note that an iteration, as used in \cite{MS95,Schn99JMPA}
with a sequence of suitable chosen $ \delta_j $s, 
is not possible
in case of periodic boundary conditions. 
}
\end{remark}

\appendix

\section{The analytic set-up}
\label{appA2}

This section contains a few preparations for the subsequent proofs 
of the attractivity and approximation result.

i) It turns out to be advantageous that all variables 
in \eqref{imag1}-\eqref{imag3}
have the same regularity, i.e., 
we introduce the new variable $ \check{v} $ by 
$$ v =  \langle \partial_x \rangle \check{v}  = (1-\partial_x^2)^{1/2} \check{v}  $$ such that \eqref{imag1}-\eqref{imag3}
becomes 
 \begin{eqnarray} \label{imag1a}
\partial_t \check{u}_1 & = & \lambda_1 \check{u}_1 + 
\check{f}_{1,n}(\check{u}_1,u_s,\check{v} )
 , \\
\partial_t u_s & = & \Lambda_s u_s +\check{f}_{s,1}(\check{u}_1,u_s,\check{v} ),\\
\partial_t \check{v}  & = & \Lambda_v \check{v}   + \partial_x^2 \check{g}_1(\check{u}_1,u_s), \label{imag3a}
\end{eqnarray}
with 
\begin{eqnarray*} 
\check{f}_{1,n}(\check{u}_1,u_s,\check{v} ) & = & \check{f}_1(\check{u}_1,u_s,\langle \partial_x \rangle \check{v} )  \\ & = & \mathcal{O}(|\check{u}_1|^3+ |\check{u}_1||u_s| +  |u_s|^2+(|\check{u}_1|+|u_s|) |\check{v} |) , \\
\check{f}_{s,n}(\check{u}_1,u_s,\check{v} ) & = & \check{f}_s(\check{u}_1,u_s,\langle \partial_x \rangle \check{v})  \\ & = &  \mathcal{O}(|\check{u}_1|^2+ |u_s|^2+(|\check{u}_1|+|u_s|) |\check{v} |),\\
\check{g}_1(\check{u}_1,u_s) & = & \langle \partial_x \rangle^{-1}
\check{g}(\check{u}_1,u_s)  \\ & = &  \mathcal{O}(|\check{u}_1|^2+ |u_s|^2).
\end{eqnarray*}
As a consequence, the  nonlinearities   
$ \check{f}_{s,n} $ and $ \partial_x^2 \check{g}_1 $ 
are smooth 
mappings from $ H^{s+1}_{l,u} $ to $ H^{s}_{l,u} $.
The mapping $ \check{f}_{1,n} $ is arbitrarily smooth due to its compact 
support in Fourier space.

ii) We introduce the scaling operator 
$$ 
(S_{\delta} u)(x) = u(\delta x) 
$$ 
and the scaled spaces $ H^{s,\delta}_{l,u} = H^{s}_{l,u} $
equipped with the norm
$$ 
\| u \|_{H^{s,\delta}_{l,u} } = \| S_{1/\delta}  u \|_{H^{s}_{l,u} }
$$

iii) Before we start with estimating the linear semigroups we define the  $ H^s_{l,u} $-norm for $ s \in (0,1) $ 
by
$$
 \| u \|_{H^s_{l,u}} =
 \| u \|_{L^2_{l,u}} +  \| \partial_x  \langle   \partial_x  \rangle^{s-1}  u \|_{L^2_{l,u}} , \qquad  (s \in (0,1)).$$
For $ n \in \N_0 $ and $  s \in (0,1) $ we set 
$$
 \| u \|_{H^{n+s}_{l,u}} = \| u \|_{H^{n}_{l,u}} +  \| \partial_x^n u \|_{H^{s}_{l,u}}.
$$ 
We need
\begin{lemma} \label{lemA4}
For $ s,r \geq 0 $ there exists a $ \sigma > 0 $, a $ C \geq 1 $ such that 
for $ 0 < \varepsilon \leq \delta \leq 1 $ and all $ t \geq 0 $ the 
following estimates hold:
\begin{eqnarray*}
\| e^{\lambda_1 t} \|_{H^s_{l,u} \to H^{s+r}_{l,u}} & \leq & C (1+t^{-r/2}) e^{C \varepsilon^2 t} , \\
\| e^{\lambda_1 t} \|_{H^{s,\delta}_{l,u} \to H^{s+r,\delta}_{l,u}} & \leq & C (1+( \delta^2 t)^{-r/2})) e^{C \varepsilon^2 t} , \\
\|  e^{\Lambda_s t}  \|_{H^s_{l,u} \to H^{s+r}_{l,u}} & \leq & C e^{-\sigma t}(1+t^{-r/2})  , \\
\|  e^{\Lambda_s t}  \|_{H^{s,\delta}_{l,u} \to H^{s+r,\delta}_{l,u}} & \leq & C  e^{-\sigma t} (1+( \delta^2 t)^{-r/2})) , \\
\| e^{\Lambda_v t}  \|_{H^s_{l,u} \to H^{s+r}_{l,u}} & \leq & C (1+t^{-r/2})  , \\
\|  e^{\Lambda_v t} \|_{H^{s,\delta}_{l,u} \to H^{s+r,\delta}_{l,u}} & \leq & C   (1+( \delta^2 t)^{-r/2})) .
\end{eqnarray*}
\end{lemma}
\noindent
{\bf Proof.}
These estimates have been established in a number of papers,
cf. \cite[\S 10]{SU2017book}.
\qed

\begin{remark}{\rm 
We refrain from recalling a complete proof of Lemma \ref{lemA4}.
It is based on estimates like
$$ 
\sup_{k \in \R} | e^{-k^2 t}  (ik)^n |\leq C t^{-n/2}
$$  
for $ n \in \mathbb{N}_0 $ and on estimates like $ \widehat{\lambda}_1(k) \leq - \alpha k^2 $ for an $ \alpha > 0 $. For real-valued $ n \geq 0 $ with $ n = n_0 + s $ with $ n_0 \in \N $ and $ s \in [0,1) $ we use 
\begin{eqnarray*}
\sup_{k \in \R} | e^{-k^2 t}  k \langle k \rangle^{s-1}  | &\leq & 
\sup_{k \in \R} | e^{-k^2 t}  |k|^s |\sup_{k \in \R} | |k|^{1-s}  \langle k \rangle^{s-1}  | \leq C t^{-s/2}.
\end{eqnarray*}  
}\end{remark}

\section{Proof of the attractivity theorem \ref{thatt}}
\label{appB}

\label{secatt}

In order to prove the attractivity result we have to show that the solution $ (u,v) $  of \eqref{rd1}-\eqref{rd2}
to a small, but otherwise arbitrary initial condition $ (u_0,v_0) \in H^{n+1}_{l,u} \times H^{n}_{l,u} $ develops in
such a way that after a certain time
 it can be written in the form stated in \eqref{glans1}-\eqref{glans2}, i.e., after that time we must be able to extract functions $ A_1 $ and $ B_0 $ which are 
functions of the long spatial variable $ X = \delta x $.
For the derivation of the amplitude system \eqref{as1introg}-\eqref{as2introg}
we make a Taylor expansion w.r.t. the small perturbation parameter $ \delta $, with 
$ 0 < \varepsilon \leq  \delta \ll 1 $, and among other things we use that 
$ \partial_x^m A_1(\delta x) = \mathcal{O}(\delta^m) $ 
and $ \partial_x^m B_0(\delta x) = \mathcal{O}(\delta^m) $.
In the end  this means that we have to prove estimates such as 
$ \partial_x^m ( (E_1 u(\cdot,t)) )  = \mathcal{O}(\delta^{m+1}) $
and   $ \partial_x^m  (E_0 v(\cdot,t))  = \mathcal{O}(\delta^{m+2}) $
for $ t > 0 $ sufficiently large
with initial conditions of \eqref{rd1}-\eqref{rd2} satisfying the estimates 
assumed in Theorem \ref{thatt}.

\begin{remark} \label{rem32}{\rm 
By looking at the Fourier representation of the 
linearized problem
we see that the $ u $-solution and the $ v $-solution are exponentially damped for all wave numbers 
except around  $ k = 0$ where in physical space the solutions are of order $ \mathcal{O}(\delta) $. By nonlinear interaction 
no other modes of order $ \mathcal{O}(\delta) $ are created.
}
\end{remark}

\subsection{The first attractivity step}

We consider \eqref{rd1}-\eqref{rd2} after applying the normal form transformation 
from Section \ref{secnft}, i.e., in the following we consider 
\eqref{imag1a}-\eqref{imag3a}.

1) We start with solutions of order $ \check{u}_1 = \mathcal{O}(\delta) $,
$ u_s = \mathcal{O}(\delta) $
 and 
$  \check{v} = \mathcal{O}(\delta^2) $.
Setting 
$ \check{u}_1 = \delta \widetilde{u}_1 $,
$ u_s = \delta \widetilde{u}_s $,
and $ \check{v} = \delta^2 \widetilde{v} $
gives
 \begin{eqnarray} \label{imag1aa}
\partial_t \widetilde{u}_1 & = & \lambda_1 \widetilde{u}_1 + 
\widetilde{f}_1(\widetilde{u}_1,\widetilde{u}_s,\widetilde{v})
 , \\
\partial_t \widetilde{u}_s & = & \Lambda_s \widetilde{u}_s +\widetilde{f}_s(\widetilde{u}_1,\widetilde{u}_s,\widetilde{v}),\\
\partial_t \widetilde{v} & = & \Lambda_v \widetilde{v}  + \partial_x^2 \widetilde{g}(\widetilde{u}_1,\widetilde{u}_s), \label{imag3aa}
\end{eqnarray}
with 
\begin{eqnarray*} 
\widetilde{f}_1(\widetilde{u}_1,\widetilde{u}_s,\widetilde{v}) & = & \mathcal{O}(\delta^2 |\widetilde{u}_1|^3+ \delta |\widetilde{u}_1||\widetilde{u}_s| + \delta |\widetilde{u}_s|^2+ \delta^2 (|\widetilde{u}_1|+|\widetilde{u}_s|) |\widetilde{v}|) , \\
\widetilde{f}_s(\widetilde{u}_1,\widetilde{u}_s,\widetilde{v}) & = & \mathcal{O}(\delta |\widetilde{u}_1|^2+ \delta |\widetilde{u}_s|^2+\delta^2 (|\widetilde{u}_1|+ |\widetilde{u}_s|) |\widetilde{v}|),\\
\widetilde{g}(\widetilde{u}_1,\widetilde{u}_s) & = & \mathcal{O}(|\widetilde{u}_1|^2+ |\widetilde{u}_s|^2).
\end{eqnarray*}
Considering the variation of constant formula 
$$
\widetilde{u}_s(t)  =  e^{\Lambda_s t} \widetilde{u}_s(0)  
+ \int_0^t 
e^{\Lambda_s(t-\tau)} \widetilde{f}_s(\widetilde{u}_1,\widetilde{u}_s,\widetilde{v})(\tau) d\tau,
$$
it is easy to see that $ u_s = \mathcal{O}(\delta^2) $ for instance for 
$ t = 1/\delta^{1/4} $ using the exponential decay $ \|e^{\Lambda_s t} \|_{H^{n+1}_{l,u} \to H^{n+1}_{l,u}} \leq C e^{- \sigma t} $ for 
a $ \sigma > 0 $ independent of $ 0 < \varepsilon \leq \delta \ll 1 $ under the assumption 
that
$ \widetilde{f}_s(\widetilde{u}_1,\widetilde{u}_s,\widetilde{v})= \mathcal{O}(\delta) $ for $ t \in [0,1/\delta^{1/4}] $. However, since there is no $ \delta^{1/4} $ in front of the 
nonlinear terms in the $ \widetilde{v} $-equation we cannot guarantee that 
$ \widetilde{v}= \mathcal{O}(1) $ for $ t = \delta^{-1/4} $. In order to guarantee
this, some extra work has to be done.
Since the argument follows the arguments of next (more complicated) step 2) we assume for a moment 
that we have proved  $ \check{u}_c = \mathcal{O}(\delta) $, $ u_s = \mathcal{O}(\delta^2) $ and $ v = \mathcal{O}(\delta^2) $ for 
$ t = 1/\delta^{1/4} $ and close the gap in the proof in Remark \ref{remneu} after completing step 2).

2) We start \eqref{imag1a}-\eqref{imag3a} again, but now for initial conditions 
$ \check{u}_c = \mathcal{O}(\delta) $,
$ u_s = \mathcal{O}(\delta^2) $,
 and 
$ v = \mathcal{O}(\delta^2) $.
Setting 
$ \check{u}_1 = \delta \widetilde{u}_1 $,
$ u_s = \delta^2 \widetilde{u}_s $
and $ v = \delta^2 \widetilde{v} $.
We find now 
\begin{eqnarray} \label{imag1b}
\partial_t \widetilde{u}_1 & = & \lambda_1 \widetilde{u}_1 + 
\widetilde{f}_1(\widetilde{u}_1,\widetilde{u}_s,\widetilde{v})
 , \\
\partial_t \widetilde{u}_s & = & \Lambda_s \widetilde{u}_s +\widetilde{f}_s(\widetilde{u}_1,\widetilde{u}_s,\widetilde{v}),\\
\partial_t \widetilde{v} & = & \Lambda_v \widetilde{v}  + \partial_x^2 \widetilde{g}(\widetilde{u}_1,\widetilde{u}_s), \label{imag3b}
\end{eqnarray}
with 
\begin{eqnarray*} 
\widetilde{f}_1(\widetilde{u}_1,\widetilde{u}_s,\widetilde{v}) & = & \mathcal{O}(\delta^2 |\widetilde{u}_1|^3+ \delta^2 |\widetilde{u}_1||\widetilde{u}_s| + \delta^3 |\widetilde{u}_s|^2+ \delta^2 (|\widetilde{u}_1|+ \delta |\widetilde{u}_s|) |\widetilde{v}|) , \\
\widetilde{f}_s(\widetilde{u}_1,\widetilde{u}_s,\widetilde{v}) & = & \mathcal{O}( |\widetilde{u}_1|^2+ \delta |\widetilde{u}_1||\widetilde{u}_s| + \delta^2 |\widetilde{u}_s|^2+\delta (|\widetilde{u}_1|+\delta |\widetilde{u}_s|) |\widetilde{v}|),\\
\widetilde{g}(\widetilde{u}_1,\widetilde{u}_s) & = & \mathcal{O}(|\widetilde{u}_1|^2
+ \delta |\widetilde{u}_1||\widetilde{u}_s| 
+\delta^2  |\widetilde{u}_s|^2).
\end{eqnarray*}

Since attractivity happens on an $ \mathcal{O}(1/\delta^2) $-time scale 
we have to control the solutions of the last system on this  long time scale.
The first equation is not a problem since in front of all nonlinear terms 
there  is a factor $ \delta^2 $. In the second equation there is linear exponential damping which allows us to control all nonlinear terms 
in this equation. 
The main difficulty to control the solutions
on the long $ \mathcal{O}(1/\delta^2) $-time scale 
is the missing $ \delta^2 $
in front of the nonlinear terms in the third equation. 
In order to get this missing $ \delta^2 $ we need that $ \widetilde{u}_1 $ is two times differentiable 
w.r.t. the long space variable $ X $. In detail, we need  that these derivatives are $  \mathcal{O}(1)$-bounded. However, this is a problem since this exactly what we are going to prove and what is not true for $ t = 0 $.

a) We proceed as follows to get rid of this problem. 
We consider the variation of constant formula 
\begin{eqnarray*}
\widetilde{u}_1(t) & = & e^{\lambda_1 t} \widetilde{u}_1(0)  + \int_0^t 
e^{\lambda_1  (t-\tau)}  \widetilde{f}_1(\widetilde{u}_1,\widetilde{u}_s,\widetilde{v})(\tau)  d\tau, \\ 
\widetilde{u}_s(t)  & = &  e^{\Lambda_s t} \widetilde{u}_s(0)  
+ \int_0^t 
e^{\Lambda_s(t-\tau)}\widetilde{f}_s(\widetilde{u}_1,\widetilde{u}_s,\widetilde{v})(\tau)   d\tau, \\
\widetilde{v}(t) & = & e^{\Lambda_v t} \widetilde{v}(0)  +  \int_0^t 
e^{\Lambda_v (t-\tau)} \partial_x^2 \widetilde{g}(\widetilde{u}_1,\widetilde{u}_s)(\tau)  d\tau.
\end{eqnarray*}
For this system we are now going to establish a priori estimates 
which in combination with the local existence and uniqueness theorem 
will guarantee the long time existence of solutions on the long 
$ \mathcal{O}(1/\delta^2) $-time scale.

We set 
$$ 
S_{c,0}(t) =  \sup_{\tau \in [0,t]} \|  \widetilde{u}_1(\tau) \|_{H^{n+1}_{l,u}} ,
\qquad
S_{s,0}(t) = \sup_{\tau \in [0,t]}  \| \widetilde{u}_s(\tau)  \|_{H^{n+1}_{l,u}}  ,
$$ 
and 
$$ 
 S_{v,0}(t) = \sup_{\tau \in [0,t]}  \| \widetilde{v}(\tau) \|_{H^{n+1}_{l,u}} .
$$
Moreover, we need the quantity
\begin{eqnarray*}
S_{c,1}(t) & = &  S_{c,0}(t) + \sup_{\tau \in [0,t]}  \tau^{1/2} \|  \partial_x \widetilde{u}_1(\tau) \|_{L^{2}_{l,u}}.
\end{eqnarray*}
Before we start,
we remark that 
all $ H^s_{l,u} $-norms for $ \widetilde{u}_1 $ are equivalent due to the compact support of $ \widetilde{u}_1 $ in Fourier space.

i) We estimate
\begin{eqnarray*}
\|\widetilde{u}_1(t) \|_{H^{n+1}_{l,u}}& \leq  &
\| e^{\lambda_1 t} \widetilde{u}_1(0) \|_{H^{n+1}_{l,u}} \\ && \qquad + \int_0^t \|
e^{\lambda_1  (t-\tau)} \|_{H^{n+1}_{l,u} \to H^{n+1}_{l,u}} \|\widetilde{f}_1(\widetilde{u}_1,\widetilde{u}_s,\widetilde{v})(\tau) \|_{H^{n+1}_{l,u}} d\tau
\\ & \leq  &
C \| \widetilde{u}_1(0)  \|_{H^{n+1}_{l,u}} + C  \int_0^t 
\|   \widetilde{f}_1(\widetilde{u}_1,\widetilde{u}_s,\widetilde{v})(\tau) \|_{H^{n+1}_{l,u}}  d\tau \\
 & \leq  & C S_{c,0}(0) + C \delta^2 t (S_{c,0}(t)^3 + S_{c,0}(t)S_{s,0}(t)+  S_{c,0}(t) S_{v,0}(t) \\ 
 && \qquad \qquad \qquad   \qquad \qquad + \delta (S_{s,0}(t)^2 + S_{s,0}(t)S_{v,0}(t))), 
 \end{eqnarray*}
 where we used the semigroup estimate 
 from Lemma \ref{lemA4}
  and the bound on $  \widetilde{f}_1 $ after \eqref{imag3b}.

 ii)  Next we find 
 \begin{eqnarray*}
 \lefteqn{
 t^{1/2}  \| \partial_x  \widetilde{u}_1(t)   \|_{L^{2}_{l,u}} }
\\&\leq&
t^{1/2}  \|  \partial_x e^{\lambda_1 t} \widetilde{u}_1(0) \|_{L^{2}_{l,u}} 
\\ && \qquad 
 +  t^{1/2}  \int_0^t \| \partial_x
e^{\lambda_1  (t-\tau)} \|_{H^{n+1}_{l,u} \to L^{2}_{l,u}}  \|\widetilde{f}_1(\widetilde{u}_1,\widetilde{u}_s,\widetilde{v})(\tau) \|_{H^{n+1}_{l,u}} d\tau
\\&\leq &
 C \|  \widetilde{u}_1(0)   \|_{H^{n+1}_{l,u}}  + t^{1/2} \int_0^t  
 (t - \tau)^{-1/2} 
  \|  \widetilde{f}_1(\widetilde{u}_1,\widetilde{u}_s,\widetilde{v}) (\tau)  \|_{H^{n+1}_{l,u}}  d\tau \\ 
& \leq & C S_{c,0}(0) + C \delta^2 t^{1/2} \int_0^t 
(t - \tau)^{-1/2} d\tau \\
&& 
\times (S_{c,0}(t)^3 + S_{c,0}(t)S_{s,0}(t)+  S_{c,0}(t) S_{v,0}(t) 
 + \delta (S_{s,0}(t)^2 + S_{s,0}(t)S_{v,0}(t)))
\\ 
 & \leq  & C S_{c,0}(0)  + C \delta^2 t (S_{c,0}(t)^3 + S_{c,0}(t)S_{s,0}(t)+  S_{c,0}(t) S_{v,0}(t) 
 \\  && \qquad \qquad \qquad  \qquad \qquad 
 + \delta (S_{s,0}(t)^2 + S_{s,0}(t)S_{v,0}(t))), 
 \end{eqnarray*}
 where we used the semigroup estimate from Lemma \ref{lemA4}
with $ r = 1 $  and again the bound on $  \widetilde{f}_1 $ after \eqref{imag3b}.

iii) For the exponentially damped part 
we use that 
$$
 \int_0^{t} \|  e^{\Lambda_s  (t-\tau)}  \|_{H^n_{l,u} \to H^{n+1}_{l,u} } d\tau
\leq 
 \int_0^{t} e^{- \sigma (t-\tau)} (1+ (t-\tau)^{-1/2})d\tau = \mathcal{O}(1) $$
 uniformly in $ t \geq 0 $,
and the  bound on $  \widetilde{f}_s $ after \eqref{imag3b},
and so we find 
\begin{eqnarray*}
\|\widetilde{u}_s(t) \|_{H^{n+1}_{l,u} }  & \leq &
 \|e^{\Lambda_s t} \widetilde{u}_s(0)   \|_{H^{n+1}_{l,u} } 
\\ && \qquad + \int_0^t 
\|e^{\Lambda_s(t-\tau)} \|_{H^n_{l,u} \to H^{n+1}_{l,u} } \|  \widetilde{f}_s(\widetilde{u}_1,\widetilde{u}_s,\widetilde{v})(\tau) \|_{H^n_{l,u} }  d\tau
\\ 
& \leq &
 \|  \widetilde{u}_s(0)   \|_{H^{n+1}_{l,u} } 
\\ && \qquad + \int_0^t 
e^{- \sigma (t-\tau)} (1+ (t-\tau)^{-1/2})  \|  \widetilde{f}_s(\widetilde{u}_1,\widetilde{u}_s,\widetilde{v})(\tau) \|_{H^n_{l,u} }  d\tau
\\
& \leq & 
C S_{s,0}(0) + C (S_{c,0}(t)^2 +  \delta S_{c,0}(t) S_{s,0}(t)+ \delta^2 S_{s,0}(t)^2 \\&& \qquad + \delta  S_{v,0}(t) (S_{c,0}(t)+ \delta S_{s,0}(t))).
\end{eqnarray*}

iv) The estimates for   the $ \widetilde{v} $-variable are obtained from 
\begin{eqnarray*}
\|\widetilde{v}(t) \|_{H^{n+1}_{l,u} }  &  \leq &  C \|\widetilde{v}(0) \|_{H^{n+1}_{l,u} } 
\\ && + \int_0^t 
\| e^{\Lambda_v (t-\tau)} \partial_x  \|_{H^{n+1}_{l,u}  \to H^{n+1}_{l,u} } \| \partial_x \widetilde{g}(\widetilde{u}_1,\widetilde{u}_s)(\tau) \|_{H^{n+1}_{l,u} } d\tau \\ 
& \leq & C S_{v,0}(0)  + C  \int_0^t(t - \tau)^{-1/2} \tau^{-1/2} d\tau  S_{c,0}(t) S_{c,1}(t) \\
&&  
 +C  \int_0^t(t - \tau)^{-1/2}  d\tau (\delta S_{c,0}(t) S_{s,0}(t) 
+ \delta^2 S_{s,0}(t) S_{s,0}(t)  )
\\
& \leq & C S_{v,0}(0)  + C S_{c,0}(t) S_{c,1}(t) 
+ C\delta t^{1/2} (S_{c,0}(t) S_{s,0}(t) 
+ \delta S_{s,0}(t)^2  ),
\end{eqnarray*}
where the semigroup term is estimated with Lemma \ref{lemA4}
with $ r = 1 $, where we used that all $ H^s_{l,u} $-norms for $ \widetilde{u}_1 $ are equivalent due to its compact support in Fourier space,
 and where we used estimates like
$$ 
|\tau^{-1/2}\tau^{1/2} \partial_x \widetilde{u}_1(\tau)| \leq \tau^{-1/2}  S_{c,1}(t).
$$ 

v) Taking the $ \sup $ w.r.t. $ t $ on the left-hand side gives the inequalities
 \begin{eqnarray*}
 S_{c,0}(t) & \leq & C S_{c,0}(0) + C T (S_{c,0}(t)^3 + S_{c,0}(t)S_{s,0}(t)+  S_{c,0}(t) S_{v,0}(t) \\
 && \qquad \qquad \qquad  + \delta (S_{s,0}(t)^2 + S_{s,0}(t)S_{v,0}(t))), \\
 S_{c,1}(t) & \leq & C S_{c,0}(0) + C T (S_{c,0}(t)^3 + S_{c,0}(t)S_{s,0}(t)+  S_{c,0}(t) S_{v,0}(t) \\ 
 && \qquad \qquad \qquad  + \delta (S_{s,0}(t)^2 + S_{s,0}(t)S_{v,0}(t))), \\
 S_{s,0}(t) & \leq & C S_{s,0}(0) + C (S_{c,0}(t)^2 +  \delta S_{c,0}(t) S_{s,0}(t)+ \delta^2 S_{s,0}(t)^2 \\ && \qquad + \delta  S_{v,0}(t) (S_{c,0}(t)+ \delta S_{s,0}(t))), \\
 S_{v,0}(t) & \leq &  C S_{v,0}(0) + C S_{c,0}(t) S_{c,1}(t) + C  T^{1/2} (S_{c,0}(t) S_{s,0}(t) + 
\delta S_{s,0}(t)^2 ) .
\end{eqnarray*}
For $ \delta > 0 $ and $ T > 0 $ sufficiently small the last two inequalities allow to estimate
$ S_{s,0}(t)  $ and $ S_{v,0}(t)  $ in terms of 
$ S_{s,0}(0)  $, $ S_{v,0}(0)  $, $ S_{c,0}(t) $, and $ S_{c,1}(t) $.
Replacing then $ S_{s,0}(t)  $ and $ S_{v,0}(t)  $ in the first two  inequalities by these estimates
and choosing $ \delta_0 > 0 $ and $ T_1 = \mathcal{O}(1) $ sufficiently small, gives the existence of a $ C_1 = 
 \mathcal{O}(1) $ with 
\begin{equation} \label{ineq1cc}
 S_{c,0}(t) + S_{c,1}(t) +S_{s,0}(t)  +  S_{v,0}(t)  \leq C_1 
\end{equation}
for all $ t \in [0,T_1/\delta^2] $ and $ \delta \in (0,\delta_0) $.
\begin{remark}\label{remneu}{\rm 
It remains to close Step 1), i.e., we have to prove that 
$ \check{u}_1 = \mathcal{O}(\delta) $,
$ u_s = \mathcal{O}(\delta) $
 and 
$  \check{v} = \mathcal{O}(\delta^2) $ tor $ t \in [0,1/\delta^{1/4}] $.
In order to do so, we follow the argument in Step 2) but now with 
$ u_s = \delta \widetilde{u}_s $ instead of $ u_s = \delta^2 \widetilde{u}_s $.
Moreover, we set 
$$ 
S_{u,0}(t) =  \sup_{\tau \in [0,t]} \|  \widetilde{u}_1(\tau) \|_{H^{n+1}_{l,u}} + \sup_{\tau \in [0,t]}  \| \widetilde{u}_s(\tau)  \|_{H^{n+1}_{l,u}}  ,
$$ 
and 
\begin{eqnarray*}
S_{u,1}(t) & = &  S_{u,0}(t) + \sup_{\tau \in [0,t]}  \tau^{1/2} \|  \partial_x \widetilde{u}_1(\tau) \|_{L^{2}_{l,u}}
+ \sup_{\tau \in [0,t]}  \tau^{1/2} \|  \partial_x \widetilde{u}_s(\tau) \|_{L^{2}_{l,u}}
\end{eqnarray*}
With exactly the same calculations as in 2) we end up with the inequalities
 \begin{eqnarray*}
 S_{u,0}(t) & \leq & C S_{u,0}(0) + C \delta t (S_{u,0}(t)^2 + \delta  S_{u,0}(t)S_{v,0}(t)), \\
 S_{u,1}(t) & \leq & C S_{u,0}(0) +  C \delta t (S_{u,0}(t)^2 + \delta  S_{u,0}(t)S_{v,0}(t)), \\
 S_{v,0}(t) & \leq &  C S_{v,0}(0) + C S_{c,0}(t) S_{c,1}(t) .
 \end{eqnarray*}
The last  inequality allows to estimate
$ S_{v,0}(t)  $ in terms of 
$ S_{v,0}(0)  $,  $ S_{u,0}(t) $, and $ S_{u,1}(t) $.
Replacing then  $ S_{v,0}(t)  $ in the first two  inequalities by this estimate
and then choosing $ \delta_0 > 0 $  sufficiently small, gives the existence of a $ C_1 = 
 \mathcal{O}(1) $ with 
\begin{equation} \label{ineq1}
 S_{u,0}(t) + S_{u,1}(t) +  S_{v,0}(t)  \leq C_1 
\end{equation}
for all $ t \in [0,1/\delta^{1/4}] $ and $ \delta \in (0,\delta_0) $.
}\end{remark}


\subsection{The second attractivity step}

Our estimates from the first attractivity
step also guarantee that 
the solutions $ \widetilde{u}_1 $,  $ \widetilde{u}_s $, and $ \widetilde{v} $
of  \eqref{imag1}-\eqref{imag3} are $ \mathcal{O}(1) $-bounded in
$ H^{1,\delta}_{l,u} $, $ H^{n+1}_{l,u}  $, and $ H^{n+1}_{l,u} $, respectively,
on time intervals of length $ \mathcal{O}(1/\delta^2) $, for instance 
considering \eqref{ineq1} for $  t \in [T_1/(2\delta^2),T_1/\delta^2] $. 

In the next step we prove that under these assumptions $ \widetilde{u}_s $ and $ \widetilde{v} $ will be in $ H^{1/2,\delta}_{l,u} $ after an 
$ \mathcal{O}(1/\delta^2) $-time scale.
Since we have the existence and uniqueness of solutions it is sufficient to establish the bounds on this long time interval.

i) We split  
$$ 
 \widetilde{f}_s(\widetilde{u}_1,\widetilde{u}_s,\widetilde{v}) =  \widetilde{f}_{s,a}(\widetilde{u}_1) +  \widetilde{f}_{s,b}(\widetilde{u}_1,\widetilde{u}_s,\widetilde{v}),
$$
with 
\begin{eqnarray*}
\widetilde{f}_{s,a}(\widetilde{u}_1)  & = &\mathcal{O}( |\widetilde{u}_1|^2) , \\
\widetilde{f}_{s,b}(\widetilde{u}_1,\widetilde{u}_s,\widetilde{v}) & = & 
\mathcal{O}(  \delta |\widetilde{u}_1||\widetilde{u}_s| + \delta^2 |\widetilde{u}_s|^2+\delta (|\widetilde{u}_1|+\delta |\widetilde{u}_s|) |\widetilde{v}|),
\end{eqnarray*}
and find for $ t \leq  T_1/\delta^2 $ with $ T_1 = \mathcal{O}(1) $ that
\begin{eqnarray*}
\lefteqn{t^{1/4}\| \partial_x  \langle   \partial_x  \rangle^{-1/2} \widetilde{u}_s(t) \|_{H^{n+1}_{l,u} } }
\\ & \leq &
t^{1/4} \| \partial_x  \langle   \partial_x  \rangle^{-1/2}  e^{\Lambda_s t} \widetilde{u}_s(0)   \|_{H^{n+1}_{l,u} } 
\\&& \qquad + t^{1/4} \int_0^t 
\| e^{\Lambda_s(t-\tau)} \|_{H^n_{l,u} \to H^{n+1}_{l,u} } \| \partial_x  \langle   \partial_x  \rangle^{-1/2}   \widetilde{f}_{s,a}(\widetilde{u}_1)(\tau) \|_{H^n_{l,u} }  d\tau
\\
\\&& \qquad + t^{1/4} \int_0^t 
\|\partial_x \langle   \partial_x  \rangle^{-1/2}  e^{\Lambda_s(t-\tau)} \|_{H^n_{l,u} \to H^{n+1}_{l,u} } \|  \widetilde{f}_{s,b}(\widetilde{u}_1,\widetilde{u}_s,\widetilde{v})(\tau) \|_{H^n_{l,u} }  d\tau
\\
& \leq & C t^{1/4} t^{-1/4} e^{- \sigma t}  \|  \widetilde{u}_s(0)   \|_{H^{n+1}_{l,u} } 
\\&& \qquad + t^{1/4} \int_0^t 
e^{- \sigma (t-\tau)} (1+ (t-\tau)^{-1/2})  \| \partial_x  \langle   \partial_x  \rangle^{-1/2}   \widetilde{f}_{s,a}(\widetilde{u}_1)(\tau)  \|_{H^n_{l,u} }  d\tau
\\
&& \qquad + t^{1/4} \int_0^t 
e^{- \sigma (t-\tau)} (1+ (t-\tau)^{-3/4})  \|  \widetilde{f}_{s,b}(\widetilde{u}_1,\widetilde{u}_s,\widetilde{v})(\tau) \|_{H^n_{l,u} }  d\tau
\\
& \leq & 
C S_{s,0}(0) + C \delta^{1/2} t^{1/4}S_{c,1}(T_1/\delta^2)^2 \\&& \qquad  +  \delta t^{1/4} (S_{c,0}(t) S_{s,0}(t)+ \delta S_{s,0}(t)^2 +  S_{v,0}(t) (S_{c,0}(t)+ \delta S_{s,0}(t)))
\end{eqnarray*}
which is $ \mathcal{O}(1) $ for $ t = T_1/\delta^2 $.

ii) Similarly, we 
find for $ t \leq  T_1/\delta^2 $ with $ T_1 = \mathcal{O}(1) $ that
\begin{eqnarray*}
\lefteqn{t^{1/4}\| \partial_x  \langle   \partial_x  \rangle^{-1/2}\widetilde{v}(t) \|_{H^{n+1}_{l,u} } } \\ &  \leq &  C t^{1/4} \|\partial_x  \langle   \partial_x  \rangle^{-1/2}e^{\Lambda_v  t} \widetilde{v}(0) \|_{H^{n+1}_{l,u} }  
\\&& + t^{1/4} \int_0^t 
\| \partial_x  \langle   \partial_x  \rangle^{-1/2} (e^{\Lambda_v (t-\tau)} \partial_x) 
\|_{H^{n+1}_{l,u} \to H^{n+1}_{l,u} } \|
  \partial_x \widetilde{g}_b(\widetilde{u}_1,\widetilde{u}_s)(\tau) \|_{H^{n+1}_{l,u} } d\tau \\ 
& \leq & C S_{v,0}(0)  + C t^{1/4}  \int_0^t(t - \tau)^{-3/4} d\tau  \delta S_{c,1}(T_1/\delta^2)^2 \\ 
&&  
 +C t^{1/4}  \int_0^t(t - \tau)^{-3/4}  d\tau (\delta S_{c,0}(t) S_{s,0}(t) 
+ \delta^2 S_{s,0}(t) S_{s,0}(t)  )
\\
& \leq & C S_{v,0}(0)  + C \delta t^{1/2} S_{c,1}(T_1/\delta^2)^2 
+ C\delta t^{1/2} (S_{c,0}(t) S_{s,0}(t) 
+ \delta S_{s,0}(t)^2  )
\end{eqnarray*}
which is $ \mathcal{O}(1) $ for $ t = T_1/\delta^2 $.

\subsection{The attractivity induction  steps}

Our estimates from the first two attractivity
steps so far  guarantee that 
the  solutions $ \widetilde{u}_1 $,  $ \widetilde{u}_s $, and $ \widetilde{v} $
of  \eqref{imag1}-\eqref{imag3} are $ \mathcal{O}(1) $-bounded in
$ H^{1,\delta}_{l,u} $, $ H^{1/2,\delta}_{l,u}  \cap H^{n+1}_{l,u}  $, and $ H^{1/2,\delta}_{l,u}  \cap H^{n+1}_{l,u} $, respectively,
on time intervals of length $ \mathcal{O}(1/\delta^2) $. 

In the next step we prove that under these assumptions $ \widetilde{u}_1 $  will be in $ H^{3/2,\delta}_{l,u} $ after an 
$ \mathcal{O}(1/\delta^2) $-time scale. After this 
we show that this implies  that  $ \widetilde{u}_s $  and 
$ \widetilde{v} $ will be in $ H^{1,\delta}_{l,u} $ after an 
$ \mathcal{O}(1/\delta^2) $-time scale. In the next step we show 
that $ \widetilde{u}_1 $  will be in $ H^{2,\delta}_{l,u} $ after an 
$ \mathcal{O}(1/\delta^2) $-time scale, etc.
We will do this by induction.
Again it is sufficient to establish the bounds.

i) In the first step we assume that 
$$ 
U_{m,c}(t) = \sup_{\tau \in [0,t]} \|  \widetilde{u}_1(\tau) \|_{H^{m,\delta}_{l,u}}, \qquad
U_{m-1/2,s}(t) = \sup_{\tau \in [0,t]} \|  \widetilde{u}_s(\tau) \|_{H^{m-1/2,\delta}_{l,u}},
$$  
and 
$$
U_{m-1/2,v}(t) = \sup_{\tau \in [0,t]} \|  \widetilde{v}(\tau) \|_{H^{m-1/2,\delta}_{l,u}}
$$
are finite and of order $ \mathcal{O}(1) $.
We 
find for $ t \leq  T_1/\delta^2 $ with $ T_1 = \mathcal{O}(1) $
that 
\begin{eqnarray*}
 \lefteqn{
 t^{1/4}  \| \partial_x  \langle   \partial_x  \rangle^{-1/2} \widetilde{u}_1(t)   \|_{H^{m,\delta}_{l,u}} }
\\&\leq&
 t^{1/4}  \|  \partial_x  \langle   \partial_x  \rangle^{-1/2} e^{\lambda_1 t} \widetilde{u}_1(0) \|_{H^{m,\delta}_{l,u}}  
 \\&&+t^{1/4} \int_0^t \|  \partial_x  \langle   \partial_x  \rangle^{-1/2}
e^{\lambda_1  (t-\tau)} \|_{H^{m-1/2,\delta}_{l,u}  \to H^{m,\delta}_{l,u}}  \|\widetilde{f}_{1}(\widetilde{u}_1,\widetilde{u}_s,\widetilde{v})(\tau) \|_{H^{m-1/2,\delta}_{l,u}}  d\tau
\\&\leq &
 C \|  \widetilde{u}_1(0)   \|_{H^{m,\delta}_{l,u}} \\&& + t^{1/4} \int_0^t  
 (t - \tau)^{-1/4} (\delta^2(t - \tau))^{-1/4} 
  \|  \widetilde{f}_1(\widetilde{u}_1,\widetilde{u}_s,\widetilde{v}) (\tau)  \|_{H^{m-1/2,\delta}_{l,u}}   d\tau \\ 
& \leq & C \|  \widetilde{u}_1(0)   \|_{H^{m,\delta}_{l,u}} + C \delta^{3/2} t^{1/4} \int_0^t 
(t - \tau)^{-1/2} d\tau \\
&& 
\times (U_{m,c}(t)^3 + U_{m,c}(t)U_{m-1/2,s}(t)+  U_{m,c}(t) U_{m-1/2,v}(t) 
\\&& \qquad \qquad \qquad  \qquad 
 + \delta U_{m-1/2,s}(t)^2 + \delta U_{m-1/2,s}(t)U_{m-1/2,v}(t))
\\
 & \leq  &  C \|  \widetilde{u}_1(0)   \|_{H^{m,\delta}_{l,u}}  + C \delta^{3/2} t^{3/4} (U_{m,c}(t)^3 + U_{m,c}(t)U_{m-1/2,s}(t)\\ && +  U_{m,c}(t) U_{m-1/2,v}(t) 
 + \delta U_{m-1/2,s}(t)^2 + \delta U_{m-1/2,s}(t)U_{m-1/2,v}(t))
 \end{eqnarray*}
%
%
%
which is $ \mathcal{O}(1) $ for $ t = T_1/\delta^2 $
and so $ \widetilde{u}_1(t) \in H^{m+1/2,\delta}_{l,u} $ for 
$ t  = \mathcal{O}(1/\delta^2) $.

ii) In the second induction step we assume that 
$  U_{m+1/2,c}(t) $, 
$ U_{m-1/2,s}(t) $, and 
$ U_{m-1/2,v}(t) $
are finite and of order $ \mathcal{O}(1) $.
 We find for $ t \leq  T_1/\delta^2 $ with $ T_1 = \mathcal{O}(1) $ that
\begin{eqnarray*}
\lefteqn{t^{1/4}\| \partial_x  \langle   \partial_x  \rangle^{-1/2} \widetilde{u}_s(t) \|_{H^{m-1/2,\delta}_{l,u} } }
\\ & \leq &
t^{1/4} \| \partial_x  \langle   \partial_x  \rangle^{-1/2}  e^{\Lambda_s t} \widetilde{u}_s(0)   \|_{H^{m-1/2,\delta}_{l,u}  } 
\\ && + t^{1/4} \int_0^t 
\|e^{\Lambda_s(t-\tau)} \|_{ H^{m-1/2,\delta}_{l,u}  \to  H^{m-1/2,\delta}_{l,u} } 
\| \partial_x  \langle   \partial_x  \rangle^{-1/2}   \widetilde{f}_{s,a}(\widetilde{u}_1)(\tau) \|_{H^{m-1/2,\delta}_{l,u} }  d\tau
\\ &&  + t^{1/4} \int_0^t 
\|\partial_x \langle   \partial_x  \rangle^{-1/2}  e^{\Lambda_s(t-\tau)} \|_{H^{m-1/2,\delta}_{l,u}  \to  H^{m-1/2,\delta}_{l,u} } \|  \widetilde{f}_{s,b}(\widetilde{u}_1,\widetilde{u}_s,\widetilde{v})(\tau) \|_{H^{m-1/2,\delta}_{l,u} }  d\tau
\\&\leq& C t^{1/4} t^{-1/4} e^{- \sigma t}  \|  \widetilde{u}_s(0)   \|_{H^{m-1/2,\delta}_{l,u} } 
\\ && + t^{1/4} \int_0^t 
e^{- \sigma (t-\tau)}  \| \partial_x  \langle   \partial_x  \rangle^{-1/2}   \widetilde{f}_{s,a}(\widetilde{u}_1)(\tau)  \|_{H^{m-1/2,\delta}_{l,u} }  d\tau
\\
&&  + t^{1/4} \int_0^t 
e^{- \sigma (t-\tau)} (t-\tau)^{-1/4}  \|  \widetilde{f}_{s,b}(\widetilde{u}_1,\widetilde{u}_s,\widetilde{v})(\tau) \|_{H^{m-1/2,\delta}_{l,u} }  d\tau
\\
& \leq & 
C \|  \widetilde{u}_s(0)   \|_{H^{m-1/2,\delta}_{l,u} } 
+ C \delta^{1/2} t^{1/4}U_{m+1/2,c}(T_1/\delta^2)^2 \\ && \qquad  +  \delta t^{1/4} (U_{m+1/2,c}(t) U_{m-1/2,s}(t)+ \delta U_{m-1/2,s}(t)^2 \\ && \qquad \qquad  +   U_{m-1/2,v}(t) (U_{m+1/2,c}(t)(t)+ \delta U_{m-1/2,s}(t)))
\end{eqnarray*}
which is $ \mathcal{O}(1) $ for $ t = T_1/\delta^2 $ and so $ \widetilde{u}_s(t) \in H^{m,\delta}_{l,u} $ for 
$ t  = \mathcal{O}(1/\delta^2) $.

iii) In the second part of the second induction step 
we assume again that 
$  U_{m+1/2,c}(t) $, 
$ U_{m-1/2,s}(t) $, and 
$ U_{m-1/2,v}(t) $
are finite and of order $ \mathcal{O}(1) $.
 We find for $ t \leq  T_1/\delta^2 $ with $ T_1 = \mathcal{O}(1) $ that
\begin{eqnarray*}
\lefteqn{t^{1/4}\|  \partial_x  \langle   \partial_x  \rangle^{-1/2}\widetilde{v}(t) \|_{H^{n+1}_{l,u} } } \\ &  \leq &  C t^{1/4} \|  \partial_x  \langle   \partial_x  \rangle^{-1/2}e^{\Lambda_v  t} \widetilde{v}(0) \|_{H^{m-1/2,\delta}_{l,u} }  
\\&& + t^{1/4} \int_0^t 
\| \partial_x\langle   \partial_x  \rangle^{-1/2} (e^{\Lambda_v  (t-\tau)} \partial_x) 
\|_{H^{m-1/2,\delta}_{l,u}  \to H^{m-1/2,\delta}_{l,u} } \|
  \partial_x \widetilde{g}_b(\widetilde{u}_1,\widetilde{u}_s)(\tau) \|_{H^{m-1/2,\delta}_{l,u} } d\tau \\
& \leq & C \| \widetilde{v}(0) \|_{H^{m-1/2,\delta}_{l,u} }   + C t^{1/4}  \int_0^t(t - \tau)^{-3/4} d\tau  \delta U_{m+1/2,c}(T_1/\delta^2)^2 \\ 
&&  
 +C t^{1/4}  \int_0^t(t - \tau)^{-3/4}  d\tau (\delta U_{m+1/2,c}(t) U_{m-1/2,s}(t) 
+ \delta^2 U_{m-1/2,s}(t)^2  )
\\
& \leq & C \| \widetilde{v}(0) \|_{H^{m-1/2,\delta}_{l,u} }  + C \delta t^{1/2} U_{m+1/2,c}(T_1/\delta^2)^2 
\\ && + C\delta t^{1/2} (U_{m+1/2,c}(t) U_{m-1/2,s}(t) 
+ \delta U_{m-1/2,s}(t)^2  )
\end{eqnarray*}
which is $ \mathcal{O}(1) $ for $ t = T_1/\delta^2 $
 and so $ \widetilde{v}(t) \in H^{m,\delta}_{l,u} $ for 
$ t  = \mathcal{O}(1/\delta^2) $.

\subsection{Attractivity of the Ginzburg-Landau manifold}
\label{appb2}

In the first step we proved that the solutions of \eqref{rdnf1}-\eqref{rdnf3} 
develop in such a way that 
for arbitrary large but fixed  $ m $ we have 
$$ 
\delta^{-1} \| c_1|_{T_1/\delta^2} \|_{H^{m,\delta}_{l,u} \cap H^{n+1}_{l,u}}
+ \delta^{-2}\| u_s|_{T_1/\delta^2} \|_{H^{m,\delta}_{l,u} \cap H^{n+1}_{l,u}}
+ \delta^{-2}\| v|_{T_1/\delta^2} \|_{H^{m,\delta}_{l,u} \cap H^{n}_{l,u}}
= \mathcal{O}(1).
$$  
Then, we set 
\begin{eqnarray*}
 c_1 & = & \psi_{1} + \delta^2 R_{1,1} ,\\
 c_{-1} & = & \psi_{-1} + \delta^2 R_{-1,1} ,\\
 u_s & = & \psi_{s} + \delta^2 R_{s,1} ,\\
 v & = & \psi_v + \delta^3 R_{v,1}.
 \end{eqnarray*}
where $ \psi_{1} $, $ \psi_{-1} $, $ \psi_s $, and $ \psi_v $ were defined  in Section \ref{secgl42}. In the following we explain how to 
choose the $ A_{\pm,m,n} $, $ A_{s,m,n} $, and $ B_{m,n} $ initially
such that in the end  the $ \delta^2 R_{\pm 1,1} $, $ \delta^2 R_{s,1} $, and 
$ \delta^3 R_{v,1} $ will become smaller and smaller.

We start with $ A_{+,1,0}|_{T=0} = \delta^{-1}  c_1|_{t = T_1/\delta^2} $,
$ A_{-,1,0}|_{T=0} = \delta^{-1} c_{-1}|_{t = T_1/\delta^2} $,
$ A_{s,0}|_{T=0} = \delta^{-2}  u_{s}|_{t = T_1/\delta^2} $, and 
$ B_{0,0}|_{T=0} = \delta^{-2}  E_0 v|_{t = T_1/\delta^2} $.
We choose the other $ A_{\pm,m,n} $, $ A_{s,m,n} $, and  $ B_{m,n} $
as in Section \ref{secgl42}. However, by this choice we cannot guarantee that
the remaining parts of the solution 
$ \delta^2 R_{1,1} $, $ \delta^2 R_{-1,1} $, $ \delta^2 R_{s,1} $,
and $ \delta^3 R_{v,1} $ are smaller than the displayed orders w.r.t. $ \delta $.

%
%
%

These estimates can be improved by the following procedure.
For $ t = T_1/\delta^2 $ we have 
\begin{eqnarray*}
c_{1}(\delta x,T_1/\delta^2) & = & (\delta A_{+,1,0}(\delta x,0) + \delta^2 A_{+,1,1}(\delta x,0) + \delta^3 A_{+,1,2}(\delta x,0) + \ldots) \\ && + (\delta^2 A_{+,2,0}(\delta x,0) + \delta^3 A_{+,2,1}(\delta x,0) + \delta^4 A_{+,2,2}(\delta x,0) + \ldots)
\\ && + (\delta^2 A_{+,0,0}(\delta x,0) + \delta^3 A_{+,0,1}(\delta x,0) + \delta^4 A_{+,0,2}(\delta x,0) + \ldots)
\\ && + (\delta^2 A_{+,-2,0}(\delta x,0) + \delta^3 A_{+,-2,1}(\delta x,0) + \delta^4 A_{+,-2,2}(\delta x,0) + \ldots) \\ &&
+ \mathcal{O}(\delta^3).
\end{eqnarray*}
We set 
$$ 
 A_{+,1,0}(\delta x,0)=  \delta^{-1}c_{1}(\delta x,T_1/\delta^2)  .
$$ 
By this choice and the construction of the improved approximation in Section 
\ref{secgl42} we obtain initial conditions for $ \delta^2 A_{+,2,0}(\delta x,0) $,
$ \delta^2 A_{+,0,0}(\delta x,0) $, and $ \delta^2 A_{+,-2,0}(\delta x,0) $.
Therefore, for a cancelation of the $ \mathcal{O}(\delta^2) $-terms 
we set 
$$
\delta^2 A_{+,1,1}(\delta x,0) = - (\delta^2 A_{+,2,0}(\delta x,0) + \delta^2 A_{+,0,0}(\delta x,0) + \delta^2 A_{+,-2,0}(\delta x,0) ).
$$
Similarly, by the  choice of $ A_{+,1,0}(\delta x,0)  $ and $ B_{0,0}(\delta x,0) $ higher order $ \mathcal{O}(\delta^{m+2}) $-terms are determined. The $ A_{+,1,m}(\delta x,0)  $ 
can then be used to adjust the initial conditions at order $ \mathcal{O}(\delta^{m+2}) $.

Next we consider the $ B $-equation. There we have  
\begin{eqnarray*}
v(\delta x,T_1/\delta^2) & = & (\delta^2 B_{0,0}(\delta x,0) + \delta^3 B_{0,1}(\delta x,0) +  \delta^4 B_{0,2}(\delta x,0)  + \ldots) \\ && + 
(\delta^4 B_{2,0}(\delta x,0) + \delta^5 B_{2,1}(\delta x,0) +  \delta^6 B_{2,2}(\delta x,0)  + \ldots) 
\\&&
+ 
(\delta^4 B_{-2,0}(\delta x,0) + \delta^5 B_{-2,1}(\delta x,0) +  \delta^6 B_{-2,2}(\delta x,0)  + \ldots) 
\\&&
+ \mathcal{O}(\delta^5).
\end{eqnarray*}
We set 
$$ 
  B_{0,0}(\delta x,0) =  \delta^{-2}v(\delta x,T_1/\delta^2)   .
$$
By this choice, the choice of $ A_{+,1,0}(\delta x,0)  $,
and the construction of the improved approximation in Section 
\ref{secgl42} we obtain initial conditions for $ \delta^4 B_{2,0}(\delta x) $,
$ \delta^4 B_{-2,0}(\delta x) $, etc..
The $ B_{0,m}(\delta x)  $ 
can then be used to adjust the initial conditions at order $ \mathcal{O}(\delta^{m+2}) $.

Finally, we come to the $ u_s $-equation.
We have 
for $ t = T_1/\delta^2 $ that 
\begin{eqnarray*}
u_{s}(\delta x,T_1/\delta^2) & = &  (\delta^2 A_{s,2,0}(\delta x,0) + \delta^3 A_{s,2,1}(\delta x,0) + \delta^4 A_{s,2,2}(\delta x,0) + \ldots)
\\&& + (\delta^2 A_{s,0,0}(\delta x,0) + \delta^3 A_{s,0,1}(\delta x,0) + \delta^4 A_{s,0,2}(\delta x,0) + \ldots)
\\&& + (\delta^2 A_{s,-2,0}(\delta x,0) + \delta^3 A_{s,-2,1}(\delta x,0) + \delta^4 A_{s,-2,2}(\delta x,0) + \ldots) \\ &&
+ \mathcal{O}(\delta^3) + \delta^2 R_{s,0}(\delta x,0).
\end{eqnarray*}
By the  choice of $ A_{+,1,0}(\delta x,0)  $ and $ B_{0,0}(\delta x,0) $ 
the $ A_{s,2,0}(\delta x,0) $, $ A_{s,0,0}(\delta x,0) $, and $ A_{s,-2,0}(\delta x,0) $ are determined.
However, in general there is a mismatch between the solution on the 
left-hand side  and the approximation terms on the right-hand side and so we need an initial correction $ \delta^2 R_{s,0}(\delta x,0) $ on the right-hand side. Since the linear semigroup $ e^{\Lambda_s t} $ decays 
with some exponential rate the variation of constant formula 
immediately yields 
$$ \delta^2 R_{s,0}(\delta x, 1/\delta^{1/4}) =  \mathcal{O}(\delta^3) .$$ 
Then we can go on and adjust the next order  initial conditions 
in the $ c_1 $- and $ v $-equation. An iteration of this procedure
finally yields the statement  of Theorem \ref{thatt}.

\section{Proof of the approximation theorem \ref{thapp}}
\label{appC}

We consider \eqref{rd1}-\eqref{rd2} after diagonalization and application of the normal form transformation 
from Section \ref{secnft}, i.e., we consider \eqref{imag1}-\eqref{imag3}.
We introduce the error functions by
\begin{equation} \label{secCansatz}
(\check{u}_{\pm 1},u_s,v) = (\delta \check{\Psi}_{\pm 1},\delta^{2} \Psi_s,\delta^{2} \Psi_v)+  
(\delta^\theta R_{\pm 1},\delta^{\theta+1} R_s,\delta^{\theta+1} R_v)
\end{equation}
where $  (\delta \check{\Psi}_{\pm 1},\delta^{2} \Psi_s,\delta^{2} \Psi_v) $ are the components of the Ginzburg-Landau approximation for 
\eqref{imag1}-\eqref{imag3}.
We look for an $ \mathcal{O}(1) $-bound for 
$$
\| R_{\pm 1} \|_{H^{1,\delta}_{l,u}}+ \| R_s \|_{H^{n+1}_{l,u}}+\| R_v \|_{H^{n}_{l,u}}.
$$
on the long $ \mathcal{O}(1/\delta^2) $-time scale.
\begin{remark}{\rm
This choice of norms allows us to use the $ \partial_x^2 $ in front 
of nonlinearity in the $ v $-equation as follows. One $ \partial_x $ is transformed 
into a $ \delta $ by using the smoothing  of the linear semigroup, i.e.,
$$ e^{d_v \partial_x^2 t} \partial_x  = \mathcal{O}(t^{-1/2}) = \delta \mathcal{O}(T^{-1/2}) $$
where $ T = \delta^2 t $. The second $ \partial_x $ is transformed 
into a $ \delta $ by using the estimate 
$$ 
\|\partial_x u \|_{H^{m,\delta}_{l,u}} \leq  C \delta \| u \|_{H^{m+1,\delta}_{l,u}}.
$$
Thus, in sum we obtain a factor $ \delta^2 $ which allows us to bound the 
solutions on the long $ \mathcal{O}(1/\delta^2) $-time scale. 
}\end{remark}

Inserting the ansatz \eqref{secCansatz}
into \eqref{imag1}-\eqref{imag3}
and applying  the variation of constant formula
gives for the error $ (R_1, R_s,R_v) $ that
\begin{eqnarray*}
R_1(t) & = & e^{\Lambda_u t}E_1 R_1|_{t=0}
+ \int_0^t    e^{\Lambda_u (t-\tau)}E_1 N_{1}(R(\tau))
 d\tau
 ,\\
R_s(t) & = & e^{\Lambda_u t}E_s R_s|_{t=0}
+ \int_0^t e^{\Lambda_u (t-\tau)}E_s  N_s(R(\tau))
 d\tau
,\\ 
R_v(t) & = & e^{\Lambda_v t} R_v|_{t=0}
+ \int_0^t (e^{\Lambda_v (t-\tau)} \partial_x)( \partial_x N_v(R(\tau)))
 d\tau,
\end{eqnarray*}
with 
\begin{eqnarray*}
\| N_{1}(R) \|_{H^{n+1}_{l,u}} & \leq &  C (\delta^2 \tilde{R}+\delta^3 \tilde{R}^2) + C_{res} \delta^2 , \\
\| N_s(R) \|_{H^{n+1}_{l,u}} & \leq &  C  \| R_1 \|_{H^{1,\delta}_{l,u}} + C (\delta \tilde{R}+\delta \tilde{R}^2) + C_{res}, \\
\| \partial_x N_v(R) \|_{H^{n}_{l,u}} & \leq & C \delta \| R_1 \|_{H^{1,\delta}_{l,u}} + C (\delta \tilde{R}+\delta^2 \tilde{R}^2)  + C_{res} \delta^2
 \end{eqnarray*}
 where $ \tilde{R}= \tilde{R}(t) $ is defined by
\begin{equation} \label{emil2}
\tilde{R}(t):=\|{R}_{1}(t)\|_{H^{1,\delta}_{l,u}}+\|R_s(t)\|_{H^{n+1}_{l,u}}+\|R_v(t)\|_{H^{n}_{l,u}},  
\end{equation}
 and where $ C_\mathrm{Res} $ stands for the $ \mathcal{O}(1) $-constants 
 coming from the residual terms.
  
 In the following $ C_\mathrm{IR} $ denotes $ \mathcal{O}(1) $-constants 
 which are obtained when integrating the 
residual terms or $ \mathcal{O}(1) $-constants  coming from the initial conditions.
 We obtain
\begin{eqnarray*}
 \|{R}_{1}(t)\|_{H^{1,\delta}_{l,u}}
	&\leq & C_\mathrm{IR}+\int_0^t C (\frac{\delta}{\sqrt{t-\tau}} + \delta^2)(\tilde{R}(\tau)+\delta\tilde{R}(\tau)^2)\,\mathrm{d}\tau,\\
\|R_s(t)\|_{H^{n+1}_{l,u}}
	&\leq &C_\mathrm{IR}+\int_0^tCe^{-\sigma(t-\tau)}\left(1+(t-\tau)^{-1/2}\right) \\ && \qquad 
	\qquad \times
	(\|{R}_{1}(\tau)\|_{H^{1,\delta}_{l,u}}
	+\delta
	(\tilde{R}(\tau)+\tilde{R}(\tau)^2)\,\mathrm{d}\tau,\\
\|R_v(t)\|_{H^{n}_{l,u}}
	&\leq &C_\mathrm{IR}+\int_0^t\frac{C\delta}{\sqrt{t-\tau}}\left(\tilde{R}(\tau)+\delta \tilde{R}(\tau)^2\right)\,\mathrm{d}\tau
\end{eqnarray*}
using Lemma \ref{lemA4}.
Next we introduce 
\begin{eqnarray*}
q_{c}(t)&=&\sup\limits_{\tau\in[0,t]} \|{R}_{1}(\tau)\|_{H^{1,\delta}_{l,u}}, \\ 
q_{s}(t)&=&\sup\limits_{\tau\in[0,t]} \|R_s(\tau)\|_{H^{n+1}_{l,u}}, \\
q_{v}(t)&=&\sup\limits_{\tau\in[0,t]} \|R_v(\tau)\|_{H^{n}_{l,u}}.
\end{eqnarray*}
We immediately obtain 
$$
q_s(t) \leq  C_\mathrm{IR} + C q(t) + C \delta ((q(t)+q_s(t))+(q(t)+q_s(t))^2), $$
where $$ q(t) = q_c(t) + q_v(t) .$$ 
For $ C \delta(1 +(q(t)+q_s(t))) \leq 1/2 $ this yields 
$
q_s(t) \leq C (q(t) + C_\mathrm{IR}) $ 
and then as a consequence
\begin{eqnarray*}
 q_{c}(t)
	&\leq & C C_\mathrm{IR}+\int_0^tC (\frac{\delta}{\sqrt{t-\tau}} + \delta^2)(q(\tau)+\delta q(\tau)^2)\,\mathrm{d}\tau,\\
q_v(t)
&\leq & C C_\mathrm{IR}+\int_0^t\frac{C\delta}{\sqrt{t-\tau}}(q(\tau)+\delta q(\tau)^2)\,\mathrm{d}\tau.
	\end{eqnarray*}
Adding  these two inequalities yields
\begin{eqnarray*}
q(t)
	& \leq   & C C_\mathrm{IR}+\int_0^tC (\frac{\delta}{\sqrt{t-\tau}} + \delta^2)(q(\tau)+\delta q(\tau)^2)\,\mathrm{d}\tau
	\\
	& \leq   & C C_\mathrm{IR}+\int_0^t 2C (\frac{\delta}{\sqrt{t-\tau}} + \delta^2)q(\tau)\,\mathrm{d}\tau
\end{eqnarray*}
if  $ \delta q(\tau) \leq 1 $. With $ T = \delta^2 t $ and $ \tilde{q}(T) = q(t) $ this can be written as 
$$ 
\tilde{q}(T) \leq C C_\mathrm{IR}+\int_0^T2 (\frac{C}{\sqrt{T-\tilde{\tau}}}+1)\tilde{q}(\tilde{\tau})\,\mathrm{d}\tilde{\tau}.
$$ 
Since this equation is independent of $ \delta $, 
Gronwall's inequality 
immediately 
yields the existence of a constant  $ M_q = \mathcal{O}(1) $ such that 
$$ 
\sup_{T\in [0,T_0]} \tilde{q}(T) =: M_q < \infty 
$$
or equivalently
$$ 
\sup_{t\in [0,T_0/\delta^2]} q(t) =M_q < \infty .
$$
Then 
$$
q_s(t) \leq M_s := C ( C_\mathrm{IR}+ M_q).
$$
Choosing $ \delta_0 > 0 $ so small that $ \delta_0 M_q \leq 1 $ 
and $ C  \delta_0 (1+ M_q+ M_s) \leq 1/2 $ 
we proved the error estimates stated in Theorem \ref{thapp}.
\qed

\section{Proof of Theorem \ref{glglob}}

\label{appD}

For $ B $, with vanishing mean value,  we have  Poincar\'{e}'s inequality 
\begin{equation} \label{FRA2}
( \int_0^L |\partial_X^{-1} B|^2  dX)^{1/2}  \leq \frac{L}{2 \pi} ( \int_0^L |B|^2  dX)^{1/2},
\end{equation}
where $ \partial_X^{-1} B $ is defined via its Fourier transform $ \widehat{B}(k)/ik $
using $ \widehat{B}(0) = 0 $. Since 
$ \textrm{Re} \int_0^L i \gamma_0 |\partial_X A|^2  dX = 0 $
and 
$ \textrm{Re} \int_0^L i \gamma_3 |A|^4  dX = 0 $
we find 
\begin{eqnarray*}
\frac12 \frac{d}{dT} \int_0^L |A|^2  dX & = & \int_0^L - |\partial_X A |^2 + |A|^2 -|A|^4 + \beta B |A|^2 dX,
\\
\frac12 \frac{d}{dT} \int_0^L |\partial_X^{-1} B|^2 dX & = &  \int_0^L - \alpha |B|^2 -  B |A|^2 dX.
\end{eqnarray*}
In case $ \beta  > 0 $ we estimate
\begin{eqnarray*}
&& \frac12 \frac{d}{dT} \int_0^L |A|^2 + \beta  |\partial_X^{-1} B|^2 dX
\\ & \leq & \int_0^L - |\partial_X A |^2 + |A|^2 -|A|^4 + \beta B |A|^2 
 - \alpha \beta |B|^2 - \beta  B |A|^2 dX
 \\ & \leq & \int_0^L  |A|^2 -|A|^4 - \alpha \beta |B|^2  dX
  \\ & \leq & \int_0^L  1 - |A|^2  -  (2 \pi)^2 L^{-2} \alpha \beta |\partial_X^{-1}B|^2  dX,
\end{eqnarray*}
where we have used \eqref{FRA2}.
Thus, we find
$$
\limsup_{T \to \infty} \int_0^L |A|^2 + \beta  |\partial_X^{-1} B|^2 dX \leq 
L \max(1,\frac{L^2}{(2 \pi)^2 \alpha }) =: C_{\infty,0}.
$$
For estimating  the higher order derivatives we keep some of the 
negative terms in the above calculations. Doing so, we also find
\begin{eqnarray} \label{eqnew5}
&& \frac12 \frac{d}{dT} \int_0^L |A|^2 + \beta  |\partial_X^{-1} B|^2 dX
\\ & \leq & \int_0^L - |\partial_X A |^2 - \alpha \beta |B|^2
+ 1 - |A|^2  -  (2 \pi)^2 L^{-2} \alpha \beta |\partial_X^{-1}B|^2  dX.
\nonumber
\end{eqnarray}
Next we compute
\begin{eqnarray*}
&& \frac12 \frac{d}{dT} \int_0^L |\partial_X A|^2 + \beta  | B|^2 dX
\\ & \leq & \int_0^L - |\partial_X^2 A |^2 + |\partial_X A|^2 
-2 |A|^2 |\partial_X A|^2 - {\rm  Re}((1+ i \gamma_3)A^2   (\partial_X \overline{A})^2 ) 
\\ && \qquad + \beta (\partial_X   \overline{A} ) \partial_X (B A) +  \beta (\partial_X   A) \partial_X (B \overline{A}  )
 - \alpha \beta |\partial_X B|^2 - \beta (\partial_X B) \partial_X |A|^2 dX
 \\ & \leq &  \int_0^L - |\partial_X^2 A |^2 + |\partial_X A|^2 
 - \alpha \beta |\partial_X B|^2 \\ && \qquad 
 - {\rm  Re}((1+ i \gamma_3)A^2   (\partial_X \overline{A})^2 )  + 2 \beta B |\partial_X A|^2 dX.
 \end{eqnarray*}
 We add $ \gamma $ times  the inequality \eqref{eqnew5} to the last inequality.
and use that we already know  that 
 $ \int_0^L |A|^2 dX $ and  $ \int_0^L   |\partial_X^{-1} B|^2 dX $ are bounded.
 On the right hand side of the new inequality for $ \gamma > 1 $ we have the negative 
 terms 
 $$ 
 - (\gamma -1)\int_0^L  |\partial_X A |^2 dX, \quad - \int_0^L \alpha \beta \gamma |B|^2dX/2 , \quad \textrm{and} \quad  - \int_0^L  |\partial_X^2 A |^2 dX 
 $$
 which we use to estimate the remaining non-negative terms on the right hand side of the new inequality.
  
 Using Young's inequality, 
 an interpolation inequality for $  \|  \partial_X A \|_{C^0_b}^2 $,
 that
 $ \int_0^L |A|^2 dX  \leq C $ and  $ \int_0^L   |\partial_X^{-1} B|^2 dX \leq C $ for a $ C > 0 $ uniformly in time, 
we estimate for every $ \delta > 0 $ that
\begin{eqnarray*}
| \int_0^L   {\rm  Re}((1+ i \gamma_3) A^2   (\partial_X \overline{A})^2 ) dX |& \leq &
 (1+|\gamma_3|) \|  A \|_{L^2}^2  \|  \partial_X A \|_{C^0_b}^2
\\& \leq &
 C \|  \partial_X A \|_{L^2}  \|  \partial_X^2 A \|_{L^2}  \\
& \leq & \frac{1}{2 \delta} C^2  \|  \partial_X A \|_{L^2}^2   + \frac{\delta}{2}  \|  \partial_X^2 A \|_{L^2}^2 
, \\
 | \int_0^L B |\partial_X A|^2 dX | & 
   =  &   | \int_0^L |(\partial_X^{-1}B) \partial_X (|\partial_X A|^2) dX | \\
 &  \leq  & 2 \| \partial_X^{-1} B \|_{L^2} \| \partial_X^2  A \|_{L^2} \| \partial_X  A \|_{C^0_b} 
   \\ & \leq &   C \| \partial_X^2  A \|_{L^2}^{3/2} 
   \| \partial_X  A \|_{L^2}^{1/2} 
     \\ & \leq &  \frac{6}{ \delta} C^{4/3}  \|  \partial_X A \|_{L^2}^2   + \frac{\delta}{2}  \|  \partial_X^2 A \|_{L^2}^2 .
  \end{eqnarray*}
Therefore,  by choosing $ \gamma $ sufficiently large, in case of periodic boundary conditions, we have established 
a-priori estimates for $ A \in H^1 $ and $ B \in L^2 $.
Since we also have local existence and uniqueness in these spaces
for  \eqref{as5}-\eqref{as6} 
global existence in $ H^1 \times L^2 $ follows, too.
The global existence for  $ A \in H^{s+1} $ and $ B \in H^s $ follows by using the smoothing properties of the diffusion semigroup.

In case $ \beta \leq  0 $ with $ 1+  \alpha^{-1} \beta > 0 $ we proceed similarly.
However, there is no cancelation and so we compute 
\begin{eqnarray*}
&& \frac12 \frac{d}{dT} \int_0^L |A|^2 + q  |\partial_X^{-1} B|^2 dX
\\ & \leq & \int_0^L - |\partial_X A |^2 + |A|^2 -|A|^4 + \beta B |A|^2 
 - \alpha q |B|^2 - q  B |A|^2 dX
 \\ & \leq & \int_0^L  |A|^2 - r |A|^4 - r \alpha \beta |B|^2  dX
\end{eqnarray*}
under the assumption that 
we can establish an estimate 
\begin{equation} \label{conditionD}
(q - \beta)  B |A|^2  \leq (1-r) \alpha q |B|^2 + (1-r) |A|^4
\end{equation}
for an $ r \in (0,1] $. If  we have established such an estimate we can proceed as above
to establish the global existence of solutions
However,
the constant $ C_{\infty,0} $ has to be modified since we no longer have $ r = 1 $.
A simple calculation shows that the required estimate \eqref{conditionD} can be established for 
$ \alpha $, $ \beta $ satisfying 
$ 1+  \alpha^{-1} \beta > 0 $ if 
$ r > 0 $ is chosen sufficiently  small and  $ q = 2 \alpha + \beta $. 
We refrain from optimizing the bound around $ \beta = 0 $.
\qed 
\begin{remark}{\rm 
In case of periodic boundary conditions the $ H^s $-space can be embedded in  $ H^s_{l,u} $. Together with the smoothing  in case of periodic boundary conditions we have established the existence of an absorbing ball for spatially periodic 
$ A \in H^{s+1}_{l,u} $ and $ B \in H^s_{l,u} $, too. }
\end{remark}
\begin{remark}{\rm 
Dropping the periodic boundary conditions
for the problem on the real line 
the global existence question remains an open problem.
}
\end{remark}

\begin{remark}{\rm 
We expect that the condition $ 1+  \alpha^{-1} \beta > 0 $ is sharp.
The reason is as follows.
In case $ \gamma_0 = \gamma_3 = 0 $
stationary solutions can be obtained by a simple integration of the conservation law 
\eqref{as6}
giving $ \alpha B = - |A|^2 + b $, where 
$ b \in \R $ is an arbitrary constant. Inserting this into \eqref{as5}
yields 
$$ 0 =  \partial_X^2 A + (1+  \alpha^{-1} \beta b)  A   - (1+  \alpha^{-1} \beta)  A |A|^2  .$$ 
Hence, the coefficient in front of the effective nonlinear terms is only negative for 
$ 1+  \alpha^{-1} \beta > 0 $. See also \cite{Hilder22}.
}
\end{remark}

%
%
%
%
%
%
%
%
%

\bibliographystyle{alpha}
\bibliography{\jobname}

%

\end{document}